\documentclass[11pt]{amsart}
\usepackage{shortcuts}
\usepackage{amssymb}
\title{Sums of Integers with Small Prime Factors}

\author{Sarvagya Jain}
\address{Department of Mathematics and Statistics, University of Turku, Turku, Finland}
\email{sarvagyajain.math@gmail.com}

\begin{document}

\begin{abstract}
    Let $P^+(n)$ denote the largest prime factor of $n$, with $P^+(1)=1$. We prove that, for every $\beta>103/(256\sqrt e)$, there is a constant $C=C(\beta)>0$ such that every sufficiently large integer $N$ can be written as $N=a+b$, where $a,b\in\Z_{>0}$ and
    \[
        P^+(a)\leq N^\beta,\qquad P^+(b)\leq(\log N)^C.
    \]
    This strengthens Balog's result, which gives two $N^\beta$-smooth summands for $\beta>4/(9\sqrt e)$. The proof combines Balog's argument with convolution estimates that are uniform in the residue class.
\end{abstract}
\maketitle
\section{Introduction}

Goldbach's conjecture states that every even integer greater than $2$ is a sum of two primes. A similar question is whether every sufficiently large integer can be written as a sum of two integers with small prime factors. Let $P^+(n)$ denote the largest prime factor of a positive integer $n$, with $P^+(1)=1$. We call a positive integer $n$ $y$-smooth if $P^+(n)\leq y$.

Balog~\cite{Balog} proved that, for every $\beta>4/(9\sqrt e)$, every sufficiently large integer $N$ is a sum of two $N^\beta$-smooth integers. If $\Psi(x;y)$ counts the $y$-smooth integers up to $x$, then $\Psi(x;x^\beta)\gg_\beta x$ for every fixed $\beta>0$.

There are far fewer integers whose prime factors are at most a fixed power of $\log N$. More precisely, for every fixed $C>1$,
\[
    \Psi(x;(\log x)^C)=x^{1-1/C+o(1)}=o(x)
    \qquad (x\to\infty).
\]
Our main result improves Balog's exponent from $4/(9\sqrt e)$ to $103/(256\sqrt e)$ and ensures that all prime factors of the second summand are at most $(\log N)^K$, for some fixed $K>0$.

\begin{theorem}\label{thm: main theorem}
    For any $\beta > 103/(256\sqrt{e})$, there exist constants $N_0=N_0(\beta)$ and $K=K(\beta)>0$ such that, for every integer $N>N_0$, there exist $a,b\in\Z_{>0}$ with $N=a+b$ and
    \begin{align*}
        P^{+}(a) \leq N^{\beta} \quad \text{and} \quad P^{+}(b) \leq (\log N)^K.
    \end{align*}
\end{theorem}
\subsection{AI Usage}
ChatGPT/Codex was used to assist with mathematical auditing, filling in details and revision, and improvements to exposition and \LaTeX{} formatting. Initial drafts of the proofs of Lemma~\ref{lemma: smooths in short intervals} and Lemmas~\ref{lemma: convolution method} and~\ref{lemma: average over theta(l, a)} were generated using ChatGPT and Gemini Pro, respectively. The author reviewed and edited the resulting material and takes full responsibility for all arguments and results.

Conversations with LLMs seem to suggest that the level of distribution in Theorem~\ref{cor: BV type result with uniformity in the range of summation} can be improved further, thereby also improving the exponent \(103/(256\sqrt{e})\) in Theorem~\ref{thm: main theorem}. An initial draft is available on Hexagon.

\subsection{Acknowledgements}
The author would like to thank his adviser, Kaisa Matom\"aki,
for introducing the author to the paper of Balog~\cite{Balog}, for many insightful suggestions and encouragement throughout this project, and for carefully reading the manuscript. The author also thanks Alexandru Pascadi for valuable email exchanges concerning~\cite{PascadiLargeSieve} and for providing Lemma~\ref{lemma: generalisation of pascadi large sieve}. The author is grateful to Sary Drappeau and Jared Lichtman for helpful clarifications regarding~\cite{Drappeau} and~\cite{Lichtman}, respectively.
The author was supported by the Research Council of Finland (grant numbers
346307, 364214, and 370133) and by the University of Turku Graduate School's Exactus fellowship
while working on this project.

\subsection{Notation}
The letter $p$, with or without subscripts, is reserved for primes; sums and products indexed by $p$ are therefore over primes. All logarithms are natural. We write $\Z_{>0}$ for the positive integers and $(a,b)$ for the greatest common divisor of $a$ and $b$. Divisors in sums and products are positive unless otherwise specified. For a summation variable, $n\sim N$ means $N<n\leq2N$. We write $\mathbf{1}_E$ for the indicator of a condition $E$, and $|\mcA|$ for the cardinality of a finite set $\mcA$.

We use Landau's notation $F=O(G)$ and Vinogradov's notation $F\ll G$ or $G\gg F$ to mean $|F|\leq CG$ for some constant $C>0$, with $G\geq0$. Subscripts indicate permitted dependence of the implied constant, as in $O_\epsilon$ or $\ll_\epsilon$, and constants may change from line to line. For positive quantities, $F\asymp G$ means $F\ll G\ll F$. In a specified limit, $F=o(G)$ means $F/G\to0$, and $F\sim G$ means $F/G\to1$. In particular, $F(x)\ll x^{o(1)}$ means $F(x)\ll_\epsilon x^\epsilon$ for every $\epsilon>0$. The bound $A\geq x^{\Omega(1)}$ means $A\geq x^c$ for some fixed constant $c>0$.

For an integer $n$ with $|n|\geq2$, let $P^+(n)$ and $P^-(n)$ denote the largest and smallest prime divisors of $|n|$, respectively, with $P^+(\pm1)=1$ and $P^-(\pm1)=\infty$. A positive integer $n$ is $y$-smooth if $P^+(n)\leq y$. For real $x\geq0$, $y\geq2$, $q\in\Z_{>0}$, and $a\in\Z$, put
\begin{align*}
    \Psi(x;y)&:=\abs{\{n\in\Z_{>0}:n\leq x,\ P^+(n)\leq y\}},\\
    \Psi_q(x;y)&:=\abs{\{n\in\Z_{>0}:n\leq x,\ P^+(n)\leq y,\ (n,q)=1\}},\\
    \Psi(x,y;a,q)&:=\abs{\{n\in\Z_{>0}:n\leq x,\ P^+(n)\leq y,\ n\equiv a\pmod q\}}.
\end{align*}
Thus $\Psi_q$ counts smooth integers coprime to $q$, while $\Psi(x,y;a,q)$ counts those in the residue class $a\pmod q$.
We use $\phi(n)$, $\mu(n)$, and $\omega(n)$ for Euler's totient function, the M\"obius function, and the number of distinct prime divisors of $n$, respectively; $\pi(x)$ counts the primes at most $x$. For an integer $k\geq1$, the divisor function $\tau_k(n)$ counts ordered factorisations $n=n_1\cdots n_k$ into positive integers, and $\tau=\tau_2$. For arithmetic functions $f,g$, their Dirichlet convolution is $(f*g)(n):=\sum_{d\mid n}f(d)g(n/d)$.

For an integer $a$ coprime to the modulus, $\overline a$ denotes its multiplicative inverse modulo that modulus, which is understood from the expression; overbars on complex coefficients denote complex conjugation.

Notation used only in Section~\ref{sec: proof of triple conv} is introduced there as needed.

\section{Outline of the Method}
We outline the proof of Theorem~\ref{thm: main theorem}. Take $y_2=(\log N)^K$, where $K$ is a sufficiently large fixed constant, and let $y_1\geq y_2$ be chosen below. We seek a representation $N=a+b$ with $P^+(a)\leq y_1$ and $P^+(b)\leq y_2$. The reduction at the start of the proof of Theorem~\ref{thm: main theorem} allows us to assume that $P^-(N)>y_2$. Every tuple in the following set gives the required representation, so it suffices to prove that the set is nonempty:
\begin{align*}
    \mcN':= \set{(s, \ell, m, n): N = s+\ell mn,\,1\leq s<N,\, P^+(s)\leq y_2,\, P^+(\ell m n)\leq y_1,\, \ell\sim L, \,m, n\sim M},
\end{align*}
where $L,M\geq N^{\Omega(1)}$ satisfy $LM^2=N/\psi(N)$, with $\psi(N)=\log\log N$. We further assume that $2L<y_1\leq M$. Let
\begin{align*}
     \mcN:=\set{(s, \ell, m, n): N = s+\ell mn,\,1\leq s<N,\, P^+(s)\leq y_2,\,\ell\sim L, \,m, n\sim M}.
\end{align*}
Let
\begin{align*}
    \mcN_1&:=\set{(s, \ell, pm, n)\in \mcN: p\text{ prime},\, p>y_1},\text{ and }\\
    \mcN_2&:=\set{(s, \ell, m, pn)\in \mcN: p\text{ prime},\, p>y_1}.
\end{align*}
By symmetry, $\abs{\mcN_1} = \abs{\mcN_2}$. Thus, as $2L<y_1$, we have the following inequality
\begin{align}\label{eqn: sketch lower bound}
    \abs{\mcN'} &\geq \abs{\mcN}-\abs{\mcN_1}-\abs{\mcN_2}\notag\\
    &= \abs{\mcN}-2\abs{\mcN_1}.
\end{align}
Additionally, as $P^+(s)\leq y_2$ and $P^-(N)>y_2$, we have that $(N, s) = 1$. Hence, $(\ell mn, N) = 1$. Thus, it is easy to see that
\begin{align}\label{eqn: sketch count N}
    \abs{\mcN} & = \sum_{\ell\sim L}\sum_{m\sim M}\sum_{\substack{P^+(s)\leq y_2\\ s\equiv N\pmod {\ell m}\\ N-2M\ell m\leq s < N-M\ell m}}1\notag\\
    & =  \sum_{\substack{\ell\sim L\\ (\ell, N) = 1}}\sum_{\substack{m\sim M\\(m, N) = 1}}\frac{1}{\phi(\ell m)}\sum_{\substack{P^+(s)\leq y_2\\ (s, \ell m) = 1\\ N-2M\ell m\leq s < N-M\ell m}}1 +O(\mcE),
\end{align}
where
\begin{align*}
    \mcE:= \sum_{\substack{q\leq 4LM\\ (q, N) = 1}}\tau(q)\sup_{0\leq z\leq N}\Biggl|\Psi(z,y_2;N,q)- \frac{1}{\phi(q)}\Psi_q(z;y_2)\Biggr|.
\end{align*}
Suppose that, for any $A\geq 1$, we have the bound
\begin{align}\label{eqn: sketch bound E}
    \mcE\ll_A \Psi(N; y_2)(\log N)^{-A}.
\end{align}
By~\eqref{eqn: sketch bound E}, the error in~\eqref{eqn: sketch count N} is $O_A(\Psi(N;y_2)(\log N)^{-A})$.
Similarly, counting each possible prime factor $p>y_1$ and using~\eqref{eqn: sketch bound E}, we obtain
\begin{align*}
    \abs{\mcN_1} &\leq \sum_{\substack{\ell\sim L\\(\ell,N)=1}}\sum_{\substack{y_1<p\leq2M\\p\nmid N}}\sum_{\substack{m\sim M/p\\(m,N)=1}}\frac{1}{\phi(\ell pm)}
    \sum_{\substack{P^+(s)\leq y_2\\(s,\ell pm)=1\\N-2M\ell pm\leq s<N-M\ell pm}}1\\
    &\quad+O_A\pare{\frac{\Psi(N;y_2)}{(\log N)^A}}.
\end{align*}
Applying Lemmas~\ref{lemma: smooths in short intervals},~\ref{lemma: convolution method}, and~\ref{lemma: average over theta(l, a)} to the two counts above gives
\begin{align*}
    \abs{\mcN}&\gg\frac{\Psi(N;y_2)}{\psi(N)},\\
    \abs{\mcN_1}&\leq\left(\sum_{y_1<p\leq2M}\frac1p+o(1)\right)\abs{\mcN}.
\end{align*}
Thus~\eqref{eqn: sketch lower bound} and Mertens' theorem yield
\begin{align*}
    \abs{\mcN'}
    &\geq\left(1-2\sum_{y_1<p\leq2M}\frac1p+o(1)\right)\abs{\mcN}\\
    &=\left(1-2\log\frac{\log M}{\log y_1}+o(1)\right)\abs{\mcN}.
\end{align*}
If $y_1\geq M^{1/\sqrt e+\epsilon_0}$ for a fixed $0<\epsilon_0<1-1/\sqrt e$, then
\[
    1-2\log\frac{\log M}{\log y_1}
    \geq1+2\log(1/\sqrt e+\epsilon_0)>0,
\]
so $\abs{\mcN'}\gg_{\epsilon_0}\abs{\mcN}\gg\Psi(N;y_2)/\psi(N)$.
Thus we would like to minimise $M$ in order to make $y_1$ as small as possible. Suppose that~\eqref{eqn: sketch bound E} holds whenever $LM\leq N^{\vartheta}$ for some $\vartheta$ such that
$$1/2<\vartheta\quad\text{and}\quad 2\vartheta-1<(1-\vartheta)/\sqrt{e}.$$
The second inequality holds when $\vartheta<(1+\sqrt e)/(2\sqrt e+1)=0.616348\dots$. We choose $M=N^{1-\vartheta}$ and $L=N^{2\vartheta-1}/\psi(N)$, so that $LM=N^\vartheta/\psi(N)\leq N^\vartheta$ for sufficiently large $N$. For $\ell\sim L$, we then have
\[
    \ell\leq2L=\frac{2N}{\psi(N)M^2}
    =\frac{2N^{2\vartheta-1}}{\psi(N)}\ll N^{2\vartheta-1}.
\]
Since $2\vartheta-1<(1-\vartheta)/\sqrt e$, every such $\ell$ is $N^{(1-\vartheta)/\sqrt e}$-smooth for sufficiently large $N$.

The following theorem allows any $1/2<\vartheta<153/256\approx0.5977$ in~\eqref{eqn: sketch bound E}, and hence gives Theorem~\ref{thm: main theorem}.
\begin{theorem}\label{cor: BV type result with uniformity in the range of summation}
    For any $\epsilon,A>0$ and integer $k\geq1$, there exist $\delta,C>0$ such that the following holds. Let $x\geq2$ and $a_1,a_2\in\Z\setminus\{0\}$ satisfy
    \[
        \abs{a_1}\leq x^{1+\delta},\qquad
        \abs{a_2}\leq x^\delta,\qquad (a_1,a_2)=1.
    \]
    Uniformly for $(\log x)^C\leq y\leq x^{1/C}$, we have
    \begin{align*}
        &\sum_{\substack{q\leq x^{153/256-\epsilon}\\(q,a_1a_2)=1}}\tau_k(q)\sup_{0\leq z\leq x}
        \Biggl|\Psi(z,y;a_1\overline{a_2},q)
        -\frac1{\phi(q)}\Psi_q(z;y)\Biggr|
        \ll_{\epsilon,A,k}\frac{\Psi(x;y)}{(\log x)^A}.
    \end{align*}
\end{theorem}
\begin{proof}
    The proof is a standard application of Proposition~\ref{propn: triple convolution with a uniformity} below. Nevertheless, we provide some details for completeness.

    Write $\Delta(z;q):=\Psi(z,y;a_1\overline{a_2},q)-\Psi_q(z;y)/\phi(q)$. We first prove the estimate without the supremum over $z$, that is, for any fixed $B>0$,
    \begin{equation}\label{eqn: weighted endpoint distribution}
        \sum_{\substack{q\leq x^{153/256-\epsilon}\\(q,a_1a_2)=1}}
        \tau_k(q)\abs{\Delta(x;q)}
        \ll_{\epsilon,B,k}\frac{\Psi(x;y)}{(\log x)^B}.
    \end{equation}
    Following the proof of~\cite[Theorem~7.1]{PascadiLvlDist}, with $q_0=1$, $y_1=y$, $y_2=x$, and the given $k$, split the characters at conductor $D=x^\eta$, where $\eta>0$ is sufficiently small. The argument in~\cite[p.~40, paragraph beginning ``For $S_{\mathrm{small}}$'']{PascadiLvlDist}, with $A=B+1$ and summation over dyadic ranges of moduli, gives, for sufficiently large $C$,
    \[
        \sum_{q\leq x^{153/256-\epsilon}}\frac{\tau_k(q)}{\phi(q)}
        \sum_{\substack{\chi\pmod q\\1<\operatorname{cond}(\chi)\leq D}}
        \Biggl|\sum_{\substack{1\leq n\leq x\\P^+(n)\leq y}}\chi(n)\Biggr|
        \ll_{\epsilon,B,k}\frac{\Psi(x;y)}{(\log x)^B}.
    \]
    This bound is uniform in the residue: after the triangle inequality, the factors $\overline{\chi(a_1)}\chi(a_2)$ disappear as they have modulus one for $(q,a_1a_2)=1$.

    For $D\leq R\leq x^{1/2}$, Harper's estimate~\cite[Proposition~2]{HarperBV} gives
    \begin{align*}
        &\sum_{q\leq R}\frac1{\phi(q)}
        \sum_{\substack{\chi\pmod q\\\operatorname{cond}(\chi)>D}}
        \Biggl|\sum_{\substack{1\leq n\leq x\\P^+(n)\leq y}}\chi(n)\Biggr|\\
        &\qquad\ll (\log x)^{7/2}\sqrt{\Psi(x;y)}
        \left(R+x^{1/2-\eta}(\log x)^2\right),
    \end{align*}
    provided $\eta$ is sufficiently small and $C$ is sufficiently large. Taking $R=x^{1/2-c}$, with $c>0$ sufficiently small and fixed, and using $\Psi(x;y)\leq x$, gives $O(x^{1-\kappa})$ for some $\kappa>0$. This bound is independent of the residue class.

    For the remaining moduli, follow the smooth-number decomposition in the proof of~\cite[Theorem~6.4]{PascadiLvlDist}, using Proposition~\ref{propn: triple convolution with a uniformity} in place of~\cite[Proposition~6.3]{PascadiLvlDist}. Applying that proposition with cutoff $D$, and summing the dyadic blocks and separation integrals as in~\cite[Section~3.7]{Drappeau}, gives an absolute power saving.

    Combining the two ranges gives, for some $\kappa>0$,
    \[
        \sum_{\substack{q\leq x^{153/256-\epsilon}\\(q,a_1a_2)=1}}\frac1{\phi(q)}
        \Biggl|\sum_{\substack{\chi\pmod q\\\operatorname{cond}(\chi)>D}}
        \overline{\chi(a_1)}\chi(a_2)
        \sum_{\substack{1\leq n\leq x\\P^+(n)\leq y}}\chi(n)\Biggr|
        =O(x^{1-\kappa}).
    \]
    For fixed $k$, multiplication by $\max_{q\leq x}\tau_k(q)=x^{o(1)}$ leaves $O(x^{1-\kappa/2})$, which is $O(\Psi(x;y)(\log x)^{-B})$ once $C$ is sufficiently large. Thus, with suitable constants depending on $\epsilon,B,k$, we have established~\eqref{eqn: weighted endpoint distribution}.

    We now introduce the supremum, following the argument in the proof of~\cite[Corollaire~1]{Drappeau}. Put $Q=x^{153/256-\epsilon}$ and $t_0=x^{1-\gamma}$, where $\gamma>0$ will be chosen sufficiently small in terms of $\epsilon,A,k$. We first consider $0\leq z<t_0$. For $(q,a_1a_2)=1$, the triangle inequality gives
    \begin{align*}
        \sup_{0\leq z<t_0}\abs{\Delta(z;q)}
        &\leq\Psi(t_0,y;a_1\overline{a_2},q)+\frac{\Psi_q(t_0;y)}{\phi(q)}\\
        &\leq\sum_{\substack{1\leq n\leq t_0\\q\mid a_2n-a_1}}1+\frac{t_0}{\phi(q)}.
    \end{align*}
    For $a_2n\neq a_1$, the divisor bound gives
    \[
        \sum_{q\mid |a_2n-a_1|}\tau_k(q)
        =\tau_{k+1}(|a_2n-a_1|)\ll_k x^{o(1)}.
    \]
    There is at most one integer $n$ with $a_2n=a_1$. Summing over $q$ therefore gives
    \begin{align*}
        &\sum_{\substack{q\leq Q\\(q,a_1a_2)=1}}\tau_k(q)
        \sup_{0\leq z<t_0}\abs{\Delta(z;q)}\\
        &\quad\leq\sum_{\substack{1\leq n\leq t_0\\a_2n\neq a_1}}\tau_{k+1}(|a_2n-a_1|)
        +\sum_{q\leq Q}\tau_k(q)+t_0\sum_{q\leq Q}\frac{\tau_k(q)}{\phi(q)}\\
        &\quad\ll_k x^{1-\gamma+o(1)}+Q(\log x)^{k-1}.
    \end{align*}
    Here we used $\sum_{q\leq Q}\tau_k(q)/\phi(q)\ll_k(\log x)^k$. For sufficiently small $\gamma$ and sufficiently large $C$, the last bound is $O_{\epsilon,A,k}(\Psi(x;y)(\log x)^{-A})$, since $\Psi(x;y)\geq x^{1-\gamma/2}$ for sufficiently large $x$.

    We now consider the remaining range $t_0\leq z\leq x$. Apply~\eqref{eqn: weighted endpoint distribution} with $\epsilon/2$ and $B$ sufficiently large in terms of $A,k$, and denote its constants by $\delta_0,C_0$. Taking $\gamma,\delta$ sufficiently small and $C$ sufficiently large ensures that, for $t_0\leq t\leq x$,
    \[
        Q\leq t^{153/256-\epsilon/2},\qquad
        \abs{a_1}\leq t^{1+\delta_0},\qquad
        \abs{a_2}\leq t^{\delta_0},\qquad
        (\log t)^{C_0}\leq y\leq t^{1/C_0}.
    \]
    Thus~\eqref{eqn: weighted endpoint distribution} applies at every such $t$. Choose a partition
    \[
        t_0<t_1<\cdots<t_J=x,\qquad J\ll(\log x)^{A+k+1},
    \]
    such that
    \[
        \Psi(t_{j+1};y)-\Psi(t_j;y)
        \leq1+\frac{\Psi(x;y)}{(\log x)^{A+k+1}}
        \qquad(0\leq j<J).
    \]
    For $t_j\leq z\leq t_{j+1}$, the number of smooth integers $n\leq z$ in the residue class $a_1\overline{a_2}\pmod q$ lies between the corresponding counts up to $t_j$ and $t_{j+1}$. Hence
    \begin{align*}
        \Delta(z;q)&\leq\Delta(t_{j+1};q)
        +\frac{\Psi_q(t_{j+1};y)-\Psi_q(z;y)}{\phi(q)},\\
        \Delta(z;q)&\geq\Delta(t_j;q)
        -\frac{\Psi_q(z;y)-\Psi_q(t_j;y)}{\phi(q)}.
    \end{align*}
    Each numerator on the right counts smooth integers coprime to $q$ in a subinterval of $(t_j,t_{j+1}]$, so it is at most $\Psi(t_{j+1};y)-\Psi(t_j;y)$. Thus
    \[
        \abs{\Delta(z;q)}
        \leq\max\{\abs{\Delta(t_j;q)},\abs{\Delta(t_{j+1};q)}\}
        +\frac{\Psi(t_{j+1};y)-\Psi(t_j;y)}{\phi(q)}.
    \]
    Summing over $q$ and using~\eqref{eqn: weighted endpoint distribution}, together with $\sum_{q\leq Q}\tau_k(q)/\phi(q)\ll_k(\log x)^k$, we obtain
    \begin{align*}
        &\sum_{\substack{q\leq Q\\(q,a_1a_2)=1}}\tau_k(q)
        \sup_{t_0\leq z\leq x}\abs{\Delta(z;q)}\\
        &\quad\leq\sum_{j=0}^J\sum_{\substack{q\leq Q\\(q,a_1a_2)=1}}
        \tau_k(q)\abs{\Delta(t_j;q)}
        +\left(1+\frac{\Psi(x;y)}{(\log x)^{A+k+1}}\right)
        \sum_{q\leq Q}\frac{\tau_k(q)}{\phi(q)}\\
        &\quad\ll_{\epsilon,A,k}
        \frac{(J+1)\Psi(x;y)}{(\log x)^B}
        +\left(1+\frac{\Psi(x;y)}{(\log x)^{A+k+1}}\right)(\log x)^k\\
        &\quad\ll_{\epsilon,A,k}
        \frac{\Psi(x;y)}{(\log x)^{A+1}}+(\log x)^k
        \ll_{\epsilon,A,k}\frac{\Psi(x;y)}{(\log x)^A}.
    \end{align*}

    Combining the two ranges completes the proof.
\end{proof}

\section{Preliminaries}
In this section, we introduce convolution estimates and also go over results about smooth numbers that we will use in this work. For the sake of completeness, we also give estimates for certain averages of Euler products that will be used in the proof of Theorem~\ref{thm: main theorem}.
\subsection{Convolution Estimates}\label{subsec: convolution estimates}

In his original proof, Balog utilised a deep convolution estimate due to Fouvry (\cite[Théorème, condition~(C.2)]{Fouvry}) to replace the congruence condition that one obtains while applying his method. In a similar vein (see~\eqref{eqn: sketch bound E} above), our approach relies on convolution estimates, specifically incorporating recent advancements that have significantly refined these estimates.

Before stating the key convolution estimate utilised in the proof of Theorem~\ref{cor: BV type result with uniformity in the range of summation}, we establish some notation.
Let
\begin{align*}
    \omega_D(k;r):=\sum_{\substack{\chi\in \mcX_D\\\operatorname{cond}(\chi)\mid r}}\chi(k), \quad\text{where}\quad \mcX_D := \set{\chi \text{ primitive}: \operatorname{cond}(\chi)\leq D}.
\end{align*}
Note that
\begin{align*}
    \abs{\omega_D(k;r)}\leq \abs{\mcX_D}\leq D^2.
\end{align*}
Let
\begin{align*}
    \mcE_D(k; r) := \mathbf{1}_{k\equiv 1 \pmod r} - \frac{\mathbf{1}_{(k, r)=1}}{\varphi(r)}\omega_D(k;r)=\frac{\mathbf{1}_{(k, r)=1}}{\varphi(r)}\sum_{\substack{\chi\pmod r\\ \operatorname{cond}(\chi)> D}}\chi(k).
\end{align*}
As discussed in the outline, we need a level of distribution result for smooth numbers in arithmetic progressions with uniformity in the residue class in the form of Theorem~\ref{cor: BV type result with uniformity in the range of summation}.
The following convolution estimate allows us to satisfactorily bound the contribution from characters with large conductors in the proof of Theorem~\ref{cor: BV type result with uniformity in the range of summation}. We use the common admissible exponent $\theta=7/32$ for exceptional eigenvalues and unramified Hecke eigenvalues; see~\eqref{eqn: defn of theta} below.
\begin{proposition}\label{propn: triple convolution with a uniformity}
    For any $\epsilon > 0$, there exists $\delta > 0$ such that the following holds.

    Let $M, N, L, R \gg 1$, $x := MNL$, $a_1, a_2 \in \mathbb{Z} \setminus \set{0}$ with $|a_1| \leq x^{1+\delta}$, $\abs{a_2}\leq x^{\delta}$, $(a_1, a_2) = 1$, and let $(\alpha_m)_{m \sim M}, (\beta_n)_{n \sim N}, (\gamma_\ell)_{\ell \sim L}$ be $1$-bounded complex sequences.
    Suppose that, for $\theta$ as defined in~\eqref{eqn: defn of theta}, we have
    \begin{equation}\label{eqn: hypothesis for triple convolution}
        \begin{aligned}
            x^{\epsilon} \leq N, \quad &NL \leq x^{2/3-5\epsilon}, \quad L \leq x^{-\epsilon}M, \quad M \leq R \leq x^{-\epsilon}NL, \quad N^2L^3 \leq x^{1-\epsilon}R, \\
            &N^{7-4\theta}L^{4-\theta} \leq x^{2-2\theta-\epsilon}R^{2+\theta}, \quad N^{7}L^{4}\abs{a_1}^{\theta} \leq x^{2-\epsilon}R^{2}.
        \end{aligned}
    \end{equation}
    Then we have
    \begin{align*}
        \sum_{\substack{r \sim R \\ (r, a_1a_2) = 1}} \left| \sum_{\substack{m \sim M \\ n \sim N \\ \ell \sim L}} \alpha_m \beta_n \gamma_{\ell} \left( \mathbf{1}_{mn\ell \equiv a_1\overline{a_2} \pmod r} - \frac{\mathbf{1}_{(mn\ell, r) = 1}}{\phi(r)} \omega_{D}(mn\ell \overline{a_1}a_2; r) \right) \right| \ll_\epsilon \frac{x (\log x)^4}{\sqrt{\min\set{x^\delta, D}}}
    \end{align*}
    for all $D\in [1, x^{\epsilon}]$.
\end{proposition}

The proof of this proposition is deferred to Section~\ref{sec: proof of triple conv}. Estimates of this form have appeared previously in the works of Drappeau (\cite[Théorème 3]{Drappeau}) and Pascadi (\cite[Theorem 4.2]{PascadiSmoothBV}, \cite[Proposition 6.3]{PascadiLvlDist}) but they either lack the required uniformity in the residue class or give a weaker level of distribution. We use the factorisation idea of Assing--Blomer--Li~\cite{Assing-Blomer-Li}, also used by Lichtman~\cite{Lichtman}. We first prove the estimate under an auxiliary coprimality condition on the coefficients (see~\eqref{eqn: coprime support for triple convolution} below), so that the residue parameter is coprime to the full spectral level. With trivial nebentypus and our Atkin--Lehner choice of scaling matrices, the factorisation then follows directly from the inverse Hecke relation, as shown in Lemma~\ref{lemma: ABL factorisation}. We remove the auxiliary condition at the end of the proof by separating the common factors with $a_1$ and applying the multiplicative large sieve.

Theorem~\ref{cor: BV type result with uniformity in the range of summation} gives an exponent of distribution of $153/256 \approx 0.5977$, compared with $6/11 \approx 0.545$ in~\cite[Théorème 1 (1.5)]{Drappeau}.

\begin{remark}
    More generally, for every fixed $\kappa\in[0,1]$ and $\epsilon>0$, the same estimate in Theorem~\ref{cor: BV type result with uniformity in the range of summation} holds for $q\leq x^{(5-\theta\kappa)/8-\epsilon}$ when $|a_1|\leq x^{\kappa+\delta}$, with $\delta>0$ sufficiently small in terms of $\kappa,\epsilon,A,k$.
\end{remark}

\subsection{Smooth Numbers}
We recall the saddle-point method. Hildebrand and Tenenbaum in \cite[Theorems 1 and 2(i)]{Hildebrand-Tenenbaum} established that for $x \geq y \geq 2$, we have uniformly:
\begin{equation}
    \Psi(x; y) = \frac{x^\alpha \zeta(\alpha, y)}{\alpha \sqrt{2\pi(1 + (\log x)/y) \log x \log y}} \left(1 + O\left(\frac{1}{\log(u+1)} + \frac{1}{\log y}\right)\right),
\end{equation}
where $u = (\log x) / \log y$. Here, the local zeta function is defined for $\Re(s) > 0$ as:
\begin{equation*}
    \zeta(s, y) := \sum_{P^+(n)\leq y} n^{-s} = \prod_{p \leq y} (1 - p^{-s})^{-1},
\end{equation*}
and the saddle point $\alpha = \alpha(x; y) > 0$ is defined implicitly by the relation
\begin{equation}\label{eqn: saddle point defn}
    \sum_{p \leq y} \frac{\log p}{p^\alpha - 1} = \log x.
\end{equation}

It follows from~\cite[Lemmas~1 and~2]{Hildebrand-Tenenbaum} that, uniformly for $\log x<y\leq x$,
\begin{equation}\label{eqn: saddle point estimate}
    \alpha(x; y) = 1 - \frac{\log(u \log(u+1))}{\log y} + O\left(\frac{1}{\log y}\right).
\end{equation}
We remark that, by definition, $\alpha(x; y)$ is a decreasing function of $x$ for any fixed $y$, and for fixed $x$ it is nondecreasing in $y$. Additionally, for $y\geq (\log x)^K$, $\alpha(x; y)\geq 1-1/K+o(1)$.

For any $\alpha \in (0, 1)$, we define
$$\phi_\alpha(q; y) := q\prod_{p\mid q,\; p\leq y}(1-1/p^\alpha).$$
We need estimates for counting smooth numbers coprime to certain integers in our analysis. We use the notation
\[
    g_m(\sigma)=\prod_{p\mid m}(1-p^{-\sigma}),\qquad
    u_y=u+\frac{\log y}{\log(u+2)}.
\]
Following~\cite{Delabreteche-Tenenbaum}, for $m>1$ put
\[
    W_m=\log p_{\omega(m)},\qquad
    \vartheta_m=\frac{W_m}{\log y},\qquad
    E_m=\frac{\vartheta_m\{u\log(u+2)\}^{\vartheta_m}}
    {1+\vartheta_m\log(u+2)},
\]
where $p_k$ is the $k$th prime; for $m=1$, set $W_m=\vartheta_m=E_m=0$.

\begin{lemma}\label{lemma: global and local estimates for coprime smooth numbers}
    Let $x\geq y\geq2$, $y>(\log x)^2$, and let $m\geq1$ satisfy $P^+(m)\leq y$ and $\omega(m)\leq\sqrt y$. With $\alpha=\alpha(x;y)$, the following estimates hold uniformly.
    \begin{enumerate}
        \item[(i)] We have
        \[
            \Psi_m(x;y)
            =g_m(\alpha)\Psi(x;y)
            \left\{1+O\left(\frac{E_m(1+E_m)}u\right)\right\}.
        \]
        \item[(ii)] For $y\leq w\leq x$, there is a positive $b=b(x,y;w,m)\asymp1$, with absolute implied constants, such that
        \begin{align*}
            \Psi_m(w;y)
            &=g_m(\alpha)\Psi(x;y)\left(\frac wx\right)^\alpha
            \left(1-\frac{\{\log(x/w)\}^2}{2(\log x)^2}\right)^{bu}\\
            &\quad\times\left\{1+O\left(\frac1{u_y}
            +(1+E_m)\left(\frac{\log(x/w)}{\log x}+\frac{E_m}u\right)\right)\right\}.
        \end{align*}
    \end{enumerate}
\end{lemma}
\begin{proof}
    Apply Théorèmes~2.1 and~2.4(ii) of the author-corrected version of~\cite{Delabreteche-Tenenbaum}.\footnote{The author-corrected manuscript may be found at \url{https://tenenb.perso.math.cnrs.fr/PPP/PropStatFriables.pdf}.} Since $y>(\log x)^2$, their parameters satisfy $\overline{u}=u$ and $\delta=0$. In Théorème~2.4(ii), take $d=x/w$, so that $x/d=w$ and $1\leq d\leq x/y$.
\end{proof}

The next lemma gives the short-interval estimate needed below.
\begin{lemma}\label{lemma: smooths in short intervals}
    Fix $B_0>0$. Suppose that
    \[
        (\log x)^{20}\leq y\leq x^{1/20},\qquad
        1\leq q\leq x^{B_0},\qquad
        z=z(x)\to\infty,\quad z\leq(\log x)^{1/10},
    \]
    where $q$ is an integer and $x\to\infty$. Uniformly in this range, with $\alpha=\alpha(x;y)$,
    \[
        \Psi_q(x+x/z;y)-\Psi_q(x;y)
        =\frac{\phi_\alpha(q;y)}q\frac{\alpha}{z}\Psi(x;y)
        \left\{1+O_{B_0}\left(\frac1z+\frac{z}{(\log x)^{1/3}}\right)\right\}.
    \]
\end{lemma}
\begin{proof}
    Let $m=\prod_{p\mid q,\,p\leq y}p$. On $y$-smooth integers, coprimality to $m$ is equivalent to coprimality to $q$, so $\Psi_m(\,\cdot\,;y)=\Psi_q(\,\cdot\,;y)$ and $g_m(\alpha)=\phi_\alpha(q;y)/q$. Since $\omega(m)\ll_{B_0}\log x$, Lemma~\ref{lemma: global and local estimates for coprime smooth numbers} applies for large $x$. Moreover, $W_m\leq2\log\log x$ and $\vartheta_m\leq1/10$. When $u<(\log x)^{1/2}$, we also have $\vartheta_m\log\{u\log(u+2)\}=o(1)$, which gives $E_m\ll\vartheta_m$. Thus the definition of $E_m$ gives
    \[
        E_m\ll_{B_0}
        \begin{cases}
            \{u\log(u+2)\}^{1/10},&u\geq(\log x)^{1/2},\\[2pt]
            \dfrac{u\log\log x}{\log x},&u<(\log x)^{1/2}.
        \end{cases}
    \]
    These estimates, together with $u_y^{-1}\ll\sqrt{\log\log x/\log x}$, imply
    \[
        \frac{E_m(1+E_m)}u+\frac1{u_y}\ll_{B_0}(\log x)^{-1/3}.
    \]
    Lemma~\ref{lemma: global and local estimates for coprime smooth numbers}(i) therefore yields
    \begin{equation}\label{eqn: global estimate for coprime smooth numbers}
        \Psi_m(x;y)
        =g_m(\alpha)\Psi(x;y)\{1+O_{B_0}((\log x)^{-1/3})\}.
    \end{equation}

    Put $\tilde{\alpha}=\alpha(x+x/z;y)$. The same bounds hold at $x+x/z$, since $\log(x+x/z)\sim\log x$, $y>\{\log(x+x/z)\}^2$, and $\vartheta_m$ is unchanged. Thus~\eqref{eqn: global estimate for coprime smooth numbers}, with $x$ replaced by $x+x/z$, gives
    \begin{equation}\label{eqn: global estimate at upper endpoint}
        \Psi_m(x+x/z;y)
        =g_m(\tilde{\alpha})\Psi(x+x/z;y)
        \{1+O_{B_0}((\log x)^{-1/3})\}.
    \end{equation}

    Next apply Lemma~\ref{lemma: global and local estimates for coprime smooth numbers}(ii) with $x$ replaced by $x+x/z$ and $w=x$. In this application, $u,u_y,E_m$ are evaluated at $x+x/z$, and the expression inside the $O$-term satisfies
    \[
        \mathcal R:=\frac1{u_y}+(1+E_m)
        \left(\frac{\log(1+1/z)}{\log(x+x/z)}+\frac{E_m}u\right)
        \ll_{B_0}(\log x)^{-1/3},
    \]
    by the preceding bounds and $\log(1+1/z)\ll1/z$. For $b=b(x+x/z,y;x,m)\asymp1$, taking logarithms also gives
    \[
        \left(1-\frac{\{\log(1+1/z)\}^2}{2\{\log(x+x/z)\}^2}\right)^{b\log(x+x/z)/\log y}
        =1+O\left(\frac1{z^2\log x\log y}\right),
    \]
    so multiplying this factor by $1+O(\mathcal R)$ gives $1+O_{B_0}((\log x)^{-1/3})$. Lemma~\ref{lemma: global and local estimates for coprime smooth numbers}(ii) therefore yields
    \begin{equation}\label{eqn: local comparison for coprime smooth numbers}
        \Psi_m(x;y)
        =\frac{g_m(\tilde{\alpha})\Psi(x+x/z;y)}{(1+1/z)^{\tilde{\alpha}}}
        \{1+O_{B_0}((\log x)^{-1/3})\}.
    \end{equation}
    Combining~\eqref{eqn: global estimate for coprime smooth numbers}--\eqref{eqn: local comparison for coprime smooth numbers} gives
    \[
        \Psi_m(x+x/z;y)-\Psi_m(x;y)
        =g_m(\alpha)\Psi(x;y)
        \left\{(1+1/z)^{\tilde{\alpha}}-1+O_{B_0}((\log x)^{-1/3})\right\}.
    \]
    Finally, the saddle-point equation~\eqref{eqn: saddle point defn} implies
    \[
        0<-\frac{\partial\alpha(x;y)}{\partial\log x}
        =\left(\sum_{p\leq y}\frac{(\log p)^2p^\alpha}{(p^\alpha-1)^2}\right)^{-1}
        \leq\frac1{(\log2)\log x}.
    \]
    Hence $\tilde{\alpha}=\alpha+O(1/(z\log x))$. Since $1/2\leq\alpha,\tilde{\alpha}\leq1$ for sufficiently large $x$,
    \[
        (1+1/z)^{\tilde{\alpha}}-1=\frac\alpha z+O\left(\frac1{z^2}\right).
    \]
    Substituting this and $g_m(\alpha)=\phi_\alpha(q;y)/q$ proves the lemma.
\end{proof}
\subsection{Averages of Euler Products}\label{subsec: averages of Euler products}
The following two lemmas follow from a standard application of the convolution method. We suppress the dependence of the Euler products on $\alpha$ in the notation.
\begin{lemma}\label{lemma: convolution method}
    Fix $B>0$. Uniformly for $X\geq3$, $\alpha\in(1/2,1)$, $1/2\leq U\leq X^B$, $(\log X)^{20}\leq y\leq X^B$, and positive integers $\ell\leq X^B$, we have
    \[
        \sum_{n\sim U}\frac{\phi_\alpha(\ell n;y)}{\phi(\ell n)}
        =U\Theta(\ell)+O_B\left(U^{3/4}+\frac{U}{(\log X)^9}\right),
    \]
    where
    $$\Theta(\ell) = \prod_{p|\ell} \frac{1 - p^{-\alpha}}{1 - p^{-1}} \prod_{p \nmid \ell} \pare{1-\frac{p^{1-\alpha}-1}{p(p-1)}}.$$
\end{lemma}
\begin{proof}
    Put
    \[
        C_y(\ell):=\prod_{p\mid\ell}\frac{1-p^{-\alpha}\mathbf{1}_{p\leq y}}{1-p^{-1}},
    \]
    and define a multiplicative function $g$, supported on squarefree integers, by
    \[
        g(p)=\begin{cases}
            0,&p\mid\ell,\\
            (1-p^{1-\alpha})/(p-1),&p\nmid\ell,\ p\leq y,\\
            1/(p-1),&p\nmid\ell,\ p>y.
        \end{cases}
    \]
    Then
    \[
        \frac{\phi_\alpha(\ell n;y)}{\phi(\ell n)}
        =C_y(\ell)\sum_{d\mid n}g(d),\qquad
        |g(p)|\leq2p^{-\alpha}\leq2p^{-1/2}.
    \]
    In particular, uniformly in all the parameters,
    \begin{equation}\label{eqn: uniform rankin bound}
        \sum_{d\geq1}\frac{|g(d)|}{d^{3/4}}
        \leq\prod_p(1+2p^{-5/4})\ll1.
    \end{equation}
    Since the factors at primes $p\leq y$ are at most $1$, for $1\leq n\leq2U$ we have
    \begin{equation}\label{eqn: uniform bound for Euler weights}
        0<\frac{\phi_\alpha(\ell n;y)}{\phi(\ell n)}
        \leq\prod_{\substack{p\mid\ell n\\p>y}}(1-p^{-1})^{-1}
        \leq\exp\left(O_B\left(\frac{\log X}{y}\right)\right)\ll_B1.
    \end{equation}
    Here $\ell n\leq2X^{2B}$; taking $n=1$ also gives $C_y(\ell)\ll_B1$. Since
    \[
        \sum_{\substack{n\sim U\\d\mid n}}1
        =\floor{\frac{2U}d}-\floor{\frac Ud}
        =\frac Ud+O(1),
    \]
    the convolution identity gives
    \begin{align*}
        \sum_{n\sim U}\frac{\phi_\alpha(\ell n;y)}{\phi(\ell n)}
        &=C_y(\ell)\sum_{d\leq2U}g(d)\left(\frac Ud+O(1)\right)\\
        &=UC_y(\ell)\prod_{p\nmid\ell}\left(1+\frac{g(p)}p\right)+O_B(U^{3/4}).
    \end{align*}
    Here~\eqref{eqn: uniform rankin bound} bounds the floor errors and the tail introduced by completing the sum, respectively, as follows:
    \begin{align*}
        \sum_{d\leq2U}|g(d)|
        &\leq(2U)^{3/4}\sum_{d\geq1}\frac{|g(d)|}{d^{3/4}}\ll U^{3/4},\\
        U\sum_{d>2U}\frac{|g(d)|}d
        &\leq U(2U)^{-1/4}\sum_{d\geq1}\frac{|g(d)|}{d^{3/4}}\ll U^{3/4}.
    \end{align*}
    Both bounds hold for every $U\geq1/2$.
    Comparing the Euler factors, and using $\Theta(\ell)\leq1$, gives
    \begin{align*}
        \left|C_y(\ell)\prod_{p\nmid\ell}\left(1+\frac{g(p)}p\right)-\Theta(\ell)\right|
        &\ll_B\sum_{\substack{p\mid\ell\\p>y}}p^{-\alpha}+\sum_{p>y}p^{-1-\alpha}\\
        &\ll_B\frac{\log X}{\sqrt y}.
    \end{align*}
    The result follows, since $y\geq(\log X)^{20}$.
\end{proof}
\begin{lemma}\label{lemma: average over theta(l, a)}
Let $\Theta(\ell)$ be as defined in Lemma~\ref{lemma: convolution method}. As $L\to\infty$, uniformly for $\alpha\in(1/2,1)$, we have
$$ \sum_{\ell \le L} \Theta(\ell) = (1+o(1))C L,$$
where $C$ is the absolutely convergent Euler product
$$ C = \prod_p \pare{1-\frac{2p^{2-\alpha}-2p-p^{1-\alpha}+1}{p^2(p-1)}}.$$
\end{lemma}
\begin{proof}
    Define the multiplicative function
    \[
        H(n):=\prod_{p\mid n}\frac{1-p^{-\alpha}}{1-p^{-1}+p^{-2}-p^{-1-\alpha}}.
    \]
    Then
    \[
        \sum_{\ell\leq L}\Theta(\ell)
        =\Theta(1)\sum_{\ell\leq L}H(\ell).
    \]
    The function $g:=\mu*H$ is supported on squarefree integers, with
    \[
        g(p)=H(p)-1=
        \frac{(1-p^{-1})(p^{-1}-p^{-\alpha})}
        {1-p^{-1}+p^{-2}-p^{-1-\alpha}}.
    \]
    In particular, $|g(p)|\ll p^{-\alpha}$. Hence, uniformly in $\alpha$,
    \[
        \sum_{d\geq1}\frac{|g(d)|}{d^{3/4}}\ll1.
    \]
    The same convolution and completion argument as in Lemma~\ref{lemma: convolution method} yields
    \[
        \sum_{\ell\leq L}H(\ell)
        =L\prod_p\left(1+\frac{g(p)}p\right)+O(L^{3/4}).
    \]
    After multiplication by $\Theta(1)$, the local factor at each prime $p$ is
    \begin{align*}
        &\left(1-\frac1p\right)\left(1-\frac{p^{1-\alpha}-1}{p(p-1)}\right)
        +\frac1p\frac{1-p^{-\alpha}}{1-p^{-1}}\\
        &\qquad=1-\frac{2p^{2-\alpha}-2p-p^{1-\alpha}+1}{p^2(p-1)}.
    \end{align*}
    These are precisely the factors of $C$. They are bounded below by their values at $\alpha=1/2$, whose product is positive and convergent, so $C\gg1$ uniformly in $\alpha$. Since $\Theta(1)\leq1$, the error $O(L^{3/4})$ is $o(CL)$, as required.
\end{proof}

\section{Proof of Theorem \ref{thm: main theorem}}
We modify Balog's argument~\cite{Balog} by transferring the exclusion of large prime factors from the convolution coefficients to weights on the moduli. This modification is inspired by~\cite{Baker-Harman,Lichtman Shifted Primes}; see the remark following the proof.
\begin{proof}[Proof of Theorem \ref{thm: main theorem} assuming Proposition \ref{propn: triple convolution with a uniformity}]

Fix $A>0$ sufficiently large, choose $K$ sufficiently large in terms of $A$ and $\beta$, and put $y:=(\log N)^K$. Our objective is to find positive integers $a,b$ with $N=a+b$, $P^+(a)\leq N^\beta$, and $P^+(b)\leq y$. We take $N$ sufficiently large that $y\leq N^\beta$.

We first reduce the assertion for general $N$ to showing the existence assertion for sufficiently large integers $N$ satisfying $P^-(N)>(\log N)^K$. Write $N=N_1N_2$, where $N_1$ is the largest $y$-smooth divisor of $N$. Then either $N_2=1$ or $P^-(N_2)>y$.

\begin{itemize}
    \item If $N_2 = 1$, then $N$ is already $y$-smooth. For any prime $p \mid N$, we have the trivial decomposition:
    \begin{equation*}
        N = \frac{N}{p} + \frac{(p-1)N}{p},
    \end{equation*}
    where both terms are $y$-smooth.

    \item If $N_2>1$, put $y_2:=(\log N_2)^K$. Since $N_2>y=(\log N)^K$, we have $N_2\to\infty$ as $N\to\infty$. Also $y_2\leq y$ and $P^-(N_2)>y_2$. Thus applying the assertion in the proposed reduction to $N_2$ with smoothness parameter $y_2$ gives positive integers $a_2,b_2$ satisfying
    \[
        N_2=a_2+b_2,\qquad P^+(a_2)\leq N_2^\beta,\qquad P^+(b_2)\leq y_2.
    \]
    Set $a=N_1a_2$ and $b=N_1b_2$. Then $N=a+b$ and
    \[
        P^+(a)\leq\max\{y,N_2^\beta\}\leq N^\beta,
        \qquad P^+(b)\leq\max\{y,y_2\}=y.
    \]
\end{itemize}

It remains to prove existence for sufficiently large $N$ when $P^{-}(N)>y$, with $y=(\log N)^K$. Choose $1/2<\vartheta<153/256$ such that $\beta>(1-\vartheta)/\sqrt e$. Decreasing $\beta$ if necessary, we may suppose that $\beta<1-\vartheta$; also $\beta>2\vartheta-1$. Put $M:=N^{1-\vartheta}$ and let
    \begin{align*}
      \mcN:=\set{(s, \ell, m, n):N-s = \ell mn,\,1\leq s<N,\, m, n\sim N^{1-\vartheta},\,\ell\sim L:=\frac{N^{2\vartheta-1}}{\psi(N)},\text{ and } P^+(s)\leq y},
    \end{align*}
    where $\psi(N)=\log\log N$.

    As $P^{-}(N)>y$ and $P^+(s)\leq y$, we see that $(s, N)= 1$. This further implies that $(\ell mn, N) = 1$. Let
    \begin{align*}
        \mcN':= \set{(s,\ell, m, n)\in \mcN: P^+(mn\ell)\leq N^\beta}.
    \end{align*}
    Since $\ell \leq N^{2\vartheta-1}\leq N^{\beta}$,
    \begin{align*}
        \mcN'= \set{(s,\ell, m, n)\in \mcN: P^+(mn)\leq N^\beta}.
    \end{align*}

    A tuple in $\mcN'$ gives the required representation by taking $a=\ell mn$ and $b=s$. Hence it suffices to show that $\abs{\mcN'}>0$. In fact, the counting argument below gives
    \begin{align}\label{eqn: N' lower bound}
        \abs{\mcN'}\geq \Psi(N; y)N^{-o(1)}.
    \end{align}

    Observe that
    \begin{align*}
        \abs{\mcN'}\geq \abs{\mcN}-\abs{\mcN_1}-\abs{\mcN_2},
    \end{align*}
    where
    \begin{align*}
        \mcN_1:=\set{(s, \ell, pm, n)\in \mcN: p\text{ prime},\,p>N^{\beta}},
    \end{align*}
    and
    \begin{align*}
         \mcN_2:=\set{(s, \ell, m, pn)\in \mcN: p\text{ prime},\,p>N^{\beta}}.
    \end{align*}
    By symmetry, we have $\abs{\mcN_1} = \abs{\mcN_2}$. Thus,
    \begin{align}\label{eqn: combinatorial ineq}
        \abs{\mcN'}\geq \abs{\mcN}-2\abs{\mcN_1}.
    \end{align}
    Notice that
    \begin{align*}
        \abs{\mcN} = \sum_{\substack{\ell\sim L\\(\ell, N) = 1}}\sum_{\substack{m\sim N^{1-\vartheta}\\ (m, N) = 1}}\sum_{\substack{s\equiv N\pmod{\ell m}\\ N-s\sim N^{1-\vartheta}\ell m\\ P^+(s)\leq y}}1.
    \end{align*}
    Fix $\epsilon>0$ with $\vartheta<153/256-\epsilon$. We apply Theorem~\ref{cor: BV type result with uniformity in the range of summation} with $x=N$, $a_1=N$, $a_2=1$, and $k=2$, taking $K$ sufficiently large, to obtain
    \begin{align*}
        \abs{\mcN} &= \sum_{\substack{\ell\sim L\\(\ell, N) = 1}}\sum_{\substack{m\sim N^{1-\vartheta}\\ (m, N) = 1}}\frac{1}{\phi(\ell m)}\sum_{\substack{(s, \ell m) = 1\\ N-s\sim N^{1-\vartheta}\ell m\\ P^+(s)\leq y}}1 + O_A\pare{\frac{\Psi(N; y)}{(\log N)^A}}.
    \end{align*}
    For $q=\ell m$, put $h=Mq$ and $x=N-2h$. Then $h\asymp N/\psi(N)$ and $z=x/h\asymp\psi(N)$, so Lemma~\ref{lemma: smooths in short intervals} applies uniformly. Write $\alpha=\alpha(N;y)$. Apply Lemma~\ref{lemma: global and local estimates for coprime smooth numbers}(i) at $x$, and part~(ii) with $x$ replaced by $N$ and $w=x$, taking its parameter $m=\prod_{p\mid q,\,p\leq y}p$ in both applications. Since $x/N=1+O(1/\psi(N))$, the error estimates in the proof of Lemma~\ref{lemma: smooths in short intervals} give
    \begin{align*}
        \frac{\phi_{\alpha(x;y)}(q;y)}q\Psi(x;y)
        &=(1+o(1))\frac{\phi_\alpha(q;y)}q\left(\frac xN\right)^\alpha\Psi(N;y)\\
        &=(1+o(1))\frac{\phi_\alpha(q;y)}q\Psi(N;y).
    \end{align*}
    The saddle-point equation~\eqref{eqn: saddle point defn} also gives, as in the proof of Lemma~\ref{lemma: smooths in short intervals},
    \[
        \alpha(x;y)-\alpha
        =O\left(\frac{\log(N/x)}{\log x}\right)
        =O\left(\frac1{\psi(N)\log N}\right).
    \]
    Together with $h/x=(1+o(1))Mq/N$, these comparisons express the main term in Lemma~\ref{lemma: smooths in short intervals} at $N$:
    \[
        \sum_{\substack{N-2Mq\leq s<N-Mq\\P^+(s)\leq y,\ (s,q)=1}}1
        =(1+o(1))\frac{\phi_\alpha(q;y)}q\frac{\alpha Mq}{N}\Psi(N;y)+O(1).
    \]
    The $O(1)$ term accounts for the endpoints. Its total contribution after division by $\phi(q)$ is
    \[
        O\pare{\sum_{q\leq4LM}\frac{\tau_2(q)}{\phi(q)}}
        =O((\log N)^2)
        =O_A\pare{\frac{\Psi(N;y)}{(\log N)^A}}.
    \]
    We next remove the restrictions $(\ell,N)=1$ and $(m,N)=1$ from the sum of these main terms. By~\eqref{eqn: uniform bound for Euler weights}, we have $\phi_\alpha(q;y)/\phi(q)\ll1$ for $q\leq N^2$. Since $P^-(N)>y$, counting multiples of each prime divisor of $N$ bounds the resulting error by
    \begin{align*}
        \frac{\alpha M\Psi(N;y)}N
        \sum_{\substack{\ell\sim L,\ m\sim M\\(\ell m,N)>1}}
        \frac{\phi_\alpha(\ell m;y)}{\phi(\ell m)}
        &\ll\frac{LM^2\Psi(N;y)}N\sum_{p\mid N}\frac1p\\
        &\ll\frac{\Psi(N;y)}{\psi(N)}\frac{\log N}{y}.
    \end{align*}
    Here $LM^2=N/\psi(N)$. Since $\log N/y\ll(\log N)^{-9}$, removing these restrictions and applying Lemma~\ref{lemma: convolution method} with $X=N$, $U=M$, and $B=2$ produces a total error of
    \[
        O\left(\frac{\Psi(N;y)}{\psi(N)}
        \left(M^{-1/4}+(\log N)^{-9}\right)\right).
    \]
    Since $C\gg1$ uniformly for $\alpha>1/2$, this error is $o(\alpha C\Psi(N;y)/\psi(N))$. Lemma~\ref{lemma: average over theta(l, a)} now gives
    \begin{align}\label{eqn: expression for N}
        \abs{\mcN}
        &=(1+o(1))\alpha\Psi(N;y)N^{1-2\vartheta}
        \sum_{\ell\sim L}\Theta(\ell)
        +O_A\pare{\frac{\Psi(N;y)}{(\log N)^A}}\notag\\
        &=(1+o(1))\frac{\alpha C\Psi(N;y)}{\psi(N)}
        +O_A\pare{\frac{\Psi(N;y)}{(\log N)^A}}.
    \end{align}

    We now obtain an upper bound for $\abs{\mcN_1}$.
    Since $2\beta>1-\vartheta$, an integer in $(M,2M]$ has at most one prime factor exceeding $N^\beta$. Thus the multiplicity of each modulus $q=\ell pm$ is at most $\tau_2(q)$, so Theorem~\ref{cor: BV type result with uniformity in the range of summation} with $k=2$ bounds the distribution error in the count of $\mcN_1$ by $O_A(\Psi(N;y)(\log N)^{-A})$. For every $N^\beta<p\leq2M$, we have $U=M/p\geq1/2$ and $\ell p,y,U\leq N^2$ for sufficiently large $N$. Lemma~\ref{lemma: convolution method} therefore applies with $X=N$, $U=M/p$, and $B=2$ throughout this range. Partial summation using $\pi(t)\ll t/\log t$, and Mertens' theorem, give
    \begin{align*}
        \sum_{N^\beta<p\leq2M}\left(\frac Mp\right)^{3/4}
        &\leq M^{3/4}\sum_{p\leq2M}p^{-3/4}
        \ll\frac M{\log M}\ll_\vartheta\frac M{\log N},\\
        \sum_{N^\beta<p\leq2M}\frac1p
        &=\log\frac{1-\vartheta}{\beta}+o(1).
    \end{align*}
    Using the short-interval estimate as in~\eqref{eqn: expression for N} and then dropping the coprimality restrictions, Lemma~\ref{lemma: convolution method} gives
    \begin{align*}
        \abs{\mcN_1}
        &\leq(1+o(1))\alpha\Psi(N;y)N^{1-2\vartheta}
        \sum_{\ell\sim L}
        \sum_{N^\beta<p\leq2M}
        \frac{\Theta(\ell p)}p\\
        &\quad+O_\vartheta\left(\frac{\Psi(N;y)}{\psi(N)\log N}\right)
        +O_A\pare{\frac{\Psi(N;y)}{(\log N)^A}}.
    \end{align*}
    The first error is $o(\alpha C\Psi(N;y)/\psi(N))$. Moreover, $\Theta(\ell p)=(1+O(N^{-\beta/2}))\Theta(\ell)$ uniformly for these primes. Applying Lemma~\ref{lemma: average over theta(l, a)} therefore yields
    \begin{align}\label{eqn: expression for N1}
        \abs{\mcN_1}
        &\leq\left(\log\frac{1-\vartheta}{\beta}+o(1)\right)
        \frac{\alpha C\Psi(N;y)}{\psi(N)}
        +O_A\pare{\frac{\Psi(N;y)}{(\log N)^A}}.
    \end{align}
    Substituting~\eqref{eqn: expression for N} and~\eqref{eqn: expression for N1} into~\eqref{eqn: combinatorial ineq} gives
    \begin{align*}
        \abs{\mcN'}
        &\geq\left(1-2\log\frac{1-\vartheta}{\beta}+o(1)\right)
        \frac{\alpha C\Psi(N;y)}{\psi(N)}\\
        &\quad-O_A\pare{\frac{\Psi(N;y)}{(\log N)^A}}.
    \end{align*}
    Since $\alpha>1/2$, $C\gg1$, $1-2\log((1-\vartheta)/\beta)>0$, and $\psi(N)=\log\log N$, it follows that
    \[
        \abs{\mcN'}\gg\frac{\Psi(N;y)}{\psi(N)}\geq\Psi(N;y)N^{-o(1)}.
    \]
    In particular, $\abs{\mcN'}>0$ for all sufficiently large $N$, completing the proof.
\end{proof}
\begin{remark}
    The key difference between Balog's approach and ours concerns the factor
\[1-2\sum_{\substack{p\mid n\\ N^{\beta}<p\leq N^{1-\vartheta}}}1\]
which is included in the convolution coefficients in Balog's argument. Here the exclusion of large prime factors is incorporated into weights on the moduli, allowing us to apply the level of distribution results for smooth numbers in arithmetic progressions.
\end{remark}

\section{Triple convolution estimate with uniformity in residue class}\label{sec: proof of triple conv}
In this section, we establish Proposition~\ref{propn: triple convolution with a uniformity} by adapting the arguments in~\cite{Drappeau,PascadiLvlDist} to account for uniformity in the residue class. We first recall the required automorphic notation and prove the preparatory estimates.

We only give a very quick overview of the notation. For more details, we refer the reader to~\cite[Section 3]{PascadiLargeSieve}.

Let $q\in \Z_{>0}$ and $\Gamma_0(q)\subseteq\text{PSL}_2(\Z)$ be the subgroup consisting of all matrices whose bottom-left entry is divisible by $q$. The quotient $\Gamma_0(q)\backslash\mathbb H$ has finitely many cusps. We use the cusp at infinity and cusps equivalent to $1/s$, where $q=rs$ and $(r,s)=1$. For these cusps, the cusp parameter $\mu(\mathfrak a)$ of~\cite[Section~3]{PascadiLargeSieve} equals $1/q$. To a cusp $\mathfrak{a} = \tau(1/s)$ for $\tau\in \Gamma_0(q)$, we associate the scaling matrix
\begin{align}\label{eqn: choice of scaling matrix}
    \sigma_\mathfrak{a}:= \tau\cdot \begin{pmatrix}
        \sqrt{r}&-\overline{s}/\sqrt{r}\\
        s\sqrt{r}&\overline{r}\sqrt{r}
    \end{pmatrix}, \text{ for some } \overline{r}, \overline{s}\in \Z\text{ with } r\overline{r}+s\overline{s} = 1.
\end{align}
Let $f$ denote a Maass cusp form invariant under the action of $\Gamma_0(q)$. For $n\in \Z\setminus\set{0}$, let $\rho_{f, \mathfrak{a}}(n)$ denote the $n^{\text{th}}$ Fourier coefficient of $f$ at cusp $\mathfrak{a}$ with normalisation as in~\cite[Section 3]{PascadiLargeSieve}. We use the identity scaling matrix at infinity and~\eqref{eqn: choice of scaling matrix} at cusps equivalent to $1/s$.

We will restrict our treatment to the exceptional spectrum, consisting of Maass cusp forms with Laplace eigenvalue $\lambda_f<1/4$. The treatment of the remaining spectrum is standard. Set $\theta_f:=\sqrt{1-4\lambda_f}$ and fix
\begin{equation}\label{eqn: defn of theta}
    \theta:=\frac{7}{32}.
\end{equation}
Kim--Sarnak~\cite[Appendix 2]{Kim-Sarnak} gives both $\theta_f\leq\theta$ and
\[
    \abs{\lambda_f(m)}\leq\tau(m)m^{\theta/2}\qquad ((m,q)=1),
\]
where $\lambda_f(m)$ denotes the normalised unramified Hecke eigenvalue. Thus $\theta$ is a common admissible exponent for the archimedean and finite-place bounds.

All spectral sums below use a basis orthonormal for the Petersson inner product, consisting of simultaneous eigenforms for the Laplacian, the unramified Hecke operators, and the reflection $f(z)\mapsto f(-\overline z)$. Thus $f(-\overline z)=\varepsilon_f f(z)$ for some $\varepsilon_f\in\{-1,1\}$. The basis is chosen separately in each space spanned by the lifts of a newform of level dividing $q$. In the exceptional spectrum, we write $\kappa_f=i\theta_f/2$, so $\lambda_f=1/4+\kappa_f^2$. In particular,
\[
    1\geq\cosh(\pi\kappa_f)=\cos(\pi\theta_f/2)
    \geq\cos(\pi\theta/2)>0,
\]
so these normalisation factors can be absorbed into the implicit constants.

The following coprime-level form of Fourier-coefficient factorisation follows directly from the inverse Hecke relation. It implements the factorisation idea in~\cite{Assing-Blomer-Li} with our choice of scaling matrix.
\begin{lemma}\label{lemma: ABL factorisation}
    Let $q=rs$ with $(r,s)=1$, and let $f$ belong to the basis above. For $a\in\Z\setminus\{0\}$ with $(a,rs)=1$, $N\geq1/2$, any complex sequence $(a_n)_n$, and $\epsilon>0$, we have
    \[
        \abs{\sum_{n\sim N}a_n\rho_{f,1/s}(an)}
        \ll_\epsilon \abs{a}^\epsilon
        \sum_{\substack{n''\mid\abs{a}\\n''\leq2N}}
        \left(\frac{\abs{a}}{n''}\right)^{\theta/2}
        \abs{\sum_{n\sim N/n''}a_{nn''}\rho_{f,1/s}(n)}.
    \]
\end{lemma}
\begin{proof}
    The matrix $\sigma_{1/s}$ in~\eqref{eqn: choice of scaling matrix} is an Atkin--Lehner matrix for the exact divisor $r$ of $q$. For trivial nebentypus, its action commutes with the unramified Hecke operators and reflection. Consequently, for $a>0$ with $(a,q)=1$, the coefficients at $1/s$ satisfy
    \[
        \lambda_f(a)\rho_{f,1/s}(n)
        =\sum_{d\mid(a,n)}\rho_{f,1/s}(an/d^2).
    \]
    M\"obius inversion gives
    \[
        \rho_{f,1/s}(an)
        =\sum_{d\mid(a,n)}\mu(d)\lambda_f(a/d)\rho_{f,1/s}(n/d),
    \]
    and hence, on interchanging sums,
    \[
        \sum_{n\sim N}a_n\rho_{f,1/s}(an)
        =\sum_{\substack{n''\mid a\\n''\leq2N}}\mu(n'')\lambda_f(a/n'')
        \sum_{n\sim N/n''}a_{nn''}\rho_{f,1/s}(n).
    \]
    The asserted bound follows from $\abs{\lambda_f(a/n'')}\ll_\epsilon a^\epsilon(a/n'')^{\theta/2}$. For $a<0$, commutation with reflection gives $\rho_{f,1/s}(-n)=\varepsilon_f\rho_{f,1/s}(n)$, so the absolute value reduces to the case $\abs{a}>0$. The unramified Hecke relation holds for oldforms in the chosen basis as well; since the argument does not divide by $\rho_{f,1/s}(1)$, it requires no nonvanishing assumption on that coefficient.
\end{proof}
We use the Fourier conventions of~\cite{PascadiLargeSieve}. Put $e(t):=e^{2\pi it}$. For $f\in L^1(\R)$, $a\in\ell^1(\Z)$, and $\phi\in L^1(\R/\Z)$, define
\begin{align*}
    \widehat f(t)&:=\int_{\R}f(u)e(-tu)\,du\qquad(t\in\R),\\
    \widehat{a}(\alpha) &:= \sum_{n\in \Z}a_ne(-n\alpha)\qquad(\alpha\in\R/\Z),\\
    \check{\phi}(n)&:=\int_{\R/\Z}\phi(\alpha)e(n\alpha)\,d\alpha\qquad(n\in\Z).
\end{align*}
For a bounded-variation complex Borel measure $\mu$ on $\R/\Z$, likewise define
\[
    \check\mu(n):=\int_{\R/\Z}e(n\alpha)\,d\mu(\alpha)\qquad(n\in\Z).
\]
Sequences specified on a summation range are extended by zero outside that range, and we write $\norm{a_n}_2^2:=\sum_n|a_n|^2$.
We will also need some notation regarding rational approximations. Given $N>0$, let $T_{N}:(\R/\Z)^2\to \R$ denote the function
\begin{align*}
    T_{N}(\alpha, \beta):=\min_{t\in \Z_{>0}}(t+N\norm{\alpha t}+N\norm{\beta t}).
\end{align*}
Here $\|u\|:=\min_{m\in\Z}|u-m|$, and we use the shorthand $T_N(\alpha):= T_{N}(\alpha, \alpha)$. Additionally, given $N\geq 1/2$ and a bounded-variation complex Borel measure $\mu$ on $\R/\Z$, we set
\begin{align*}
    \mcI_{N}(\mu):= \iint_{(\R/\Z)^2}T_{N}(\alpha, \beta)d\abs{\mu}(\alpha)d\abs{\mu}(\beta).
\end{align*}
The following lemma is a generalisation of the spectral large sieve~\cite[Theorem 5.2]{PascadiLargeSieve} appropriate for our purposes.
\begin{lemma}\label{lemma: generalisation of pascadi large sieve}
Let $\epsilon>0$, $X, A>0$, $N\geq 1/2$, and $a, b, q\in \Z_{>0}$. Let $(a_n)_n$ be a complex sequence such that
$$a_n = F\pare{\frac{n}{N}}\check{\mu}(bn)\qquad(n\sim N)$$
for a smooth function $F:(0,4)\to\C$ satisfying $F^{(j)}\ll_j1$ for $j\geq0$ and $\mu$ a bounded-variation complex Borel measure on $\R/\Z$.
For a cusp $\mathfrak{a}$ of $\Gamma_0(q)\backslash \mathbb{H}$ with $\mu(\mathfrak{a}) = 1/q$ and the notation as above, we have
    $$\sum_{\lambda_f < 1/4} X^{\theta_f} \left| \sum_{n \sim N} a_n \rho_{f ,\mathfrak{a}}(an) \right|^2 \ll_\epsilon \bigl(qaN\max\{1,X\}\bigr)^{2\epsilon} \left(1 + \frac{aN}{q}\right) A^2,$$
provided
$$A \gg \norm{a_{n}}_2 + \frac{\sqrt{(a,q)} N}{\sqrt{q + aN}} |\mu|(\mathbb{R}/\mathbb{Z}),$$
$$X \ll \max \left(1, \frac{q}{aN}\right) \max \left(1, \frac{A^2}{\mathcal{I}_{bN}(\mu)}\right).$$
\end{lemma}
\begin{proof}
    For $0<X<1$, use $X^{\theta_f}\leq1$ and the ordinary spectral large sieve~\cite[Lemma~3.9]{PascadiLargeSieve} at length $aN$, with coefficients supported on multiples of $a$ and norm $\norm{a_n}_2\ll A$. We may therefore assume $X\geq1$.

    Let $D_b:\R/\Z\to\R/\Z$ be given by $D_b(\alpha)=b\alpha$, and let $\nu=(D_b)_*\mu$ be the pushforward measure. Then
    \[
        \check\nu(n)=\int_{\R/\Z}e(nb\alpha)\,d\mu(\alpha)=\check\mu(bn),
        \qquad |\nu|\leq(D_b)_*|\mu|.
    \]
    In particular, $|\nu|(\R/\Z)\leq|\mu|(\R/\Z)$. Also, $\|bu\|\leq b\|u\|$ gives
    \[
        T_N(b\alpha,b\beta)\leq T_{bN}(\alpha,\beta),
    \]
    and hence
    \[
        \mcI_N(\nu)
        \leq\iint_{(\R/\Z)^2}T_N(b\alpha,b\beta)\,d|\mu|(\alpha)\,d|\mu|(\beta)
        \leq\mcI_{bN}(\mu).
    \]
    Apply~\cite[Theorem 5.2]{PascadiLargeSieve} with $\nu$ in place of $\mu$, recalling that our $\theta_f=\sqrt{1-4\lambda_f}$ is twice the spectral parameter used there.
\end{proof}
The following assumption generalises \cite[Assumption 3.2]{PascadiLvlDist} and allows us to state the results succinctly. Its formulation is tailored to match the structure that arises after applying Lemma~\ref{lemma: ABL factorisation}.
\begin{assumption}\label{assumption: assumption on dispersion weights}
We say that a tuple $(q,n'',N,Z,(a_n)_{n\sim N},A_{n''},Y_{n''})$, with $q,n''\in\Z_{>0}$, $N\geq1/2$, $n''\leq2N$, $Z\gg1$, $A_{n''}\gg\norm{a_{nn''}}_2$, and $Y_{n''}>0$, satisfies this assumption if, for every $\epsilon>0$, $\xi\in\R$, and cusp $\mathfrak a$ of $\Gamma_0(q)\backslash\mathbb H$ with $\mu(\mathfrak a)=1/q$, one has
\begin{align*}
\sum_{\lambda_f<1/4}X^{\theta_f}\abs{\sum_{n\sim N/n''}a_{nn''}e\pare{\frac{nn''}{N}\xi}\rho_{f, \mathfrak{a}}(n)}^2\ll_\epsilon \pare{qNZ}^\epsilon \pare{1+\frac{N}{n''q}}A_{n''}^2
\end{align*}
for all
\begin{align*}
    X\ll \max\pare{1, \frac{qn''}{N}}\frac{Y_{n''}}{1+\abs{\xi}^2}.
\end{align*}
\end{assumption}
The following lemma verifies Assumption~\ref{assumption: assumption on dispersion weights} for arbitrary sequences and for the additive convolutions used below.
\begin{lemma}\label{lemma: assumption for dispersion weights}
Let $N\geq1/2$ and $q,n''\in\Z_{>0}$ with $n''\leq2N$.
\begin{enumerate}
\item[(i)] For any complex sequence $(a_n)_{n\sim N}$, the tuple
\[
    (q,n'',N,1,(a_n)_{n\sim N},\norm{a_{nn''}}_2,1)
\]
satisfies Assumption~\ref{assumption: assumption on dispersion weights}.
\item[(ii)] Let $L,H\gg1$, $T_0\geq1$, $\alpha_1,\alpha_2\in\R/\Z$, and $\ell_1,\ell_2\in\Z_{>0}$ satisfy $q\gg L^2$, $\ell_1,\ell_2\asymp L$, and $(\ell_1,\ell_2)=1$. For an integer $t$ with $1\leq|t|\leq T_0$, define
\[
    a_n:=\sum_{\substack{h_1,h_2\in\Z\\t(h_1\ell_1-h_2\ell_2)=n}}
    \Phi_1\pare{\frac{h_1}{H}}\Phi_2\pare{\frac{h_2}{H}}
    e(h_1\alpha_1+h_2\alpha_2)\qquad(n\sim N),
\]
where $\Phi_i:\R\to\C$ are smooth functions supported in $|u|\ll1$, with $\Phi_i^{(j)}\ll_j1$ for $j\geq0$. Then the tuple $(q,n'',N,H,(a_n)_{n\sim N},A_{n''},Y_{n''})$ satisfies Assumption~\ref{assumption: assumption on dispersion weights}, where
\begin{align*}
    A_{n''}&:=\norm{a_{nn''}}_2+\sqrt{\frac{T_0N}{n''}}
    \pare{\sqrt{\frac{H}{L}}+\frac{H}{L}},\\
    Y_{n''}&:=\max\set{1,\frac{NH}{n''L(H+L)\min_iT_H(\alpha_i)}}.
\end{align*}
\end{enumerate}
\end{lemma}
\begin{proof}
    \emph{(i)} Apply~\cite[Proposition~3.3]{PascadiLvlDist} at length $N/n''$ to the sequence $(a_{nn''})_{n\sim N/n''}$. The twist $e(n\xi/(N/n''))$ is exactly $e(nn''\xi/N)$, and $(qN/n'')^\epsilon\leq(qN)^\epsilon$, so its conclusion gives Assumption~\ref{assumption: assumption on dispersion weights} with the stated parameters.

    \emph{(ii)} Interchanging the two sets of indices and weights reduces to $t>0$. Put
    \[
        d=\frac{t}{(t,n'')},\qquad b=\frac{n''}{(t,n'')}.
    \]
    Since $(b,d)=1$, the condition $t\mid nn''$ is equivalent to $d\mid n$. Thus, for every $\xi\in\R$,
    \begin{equation}\label{eqn: dilated dispersion coefficients}
        \sum_{n\sim N/n''}a_{nn''}e\pare{\frac{nn''\xi}{N}}\rho_{f,\mathfrak a}(n)
        =\sum_{m\sim N/(dn'')}a_{dmn''}e\pare{\frac{dmn''\xi}{N}}\rho_{f,\mathfrak a}(dm).
    \end{equation}
    The sum is empty if $N/(dn'')<1/2$; otherwise we may assume $N/t\ll HL$ by the support of $\Phi_1,\Phi_2$.

    Define the measure
    \[
        d\mu(\alpha):=\prod_{i=1}^2\left(\sum_{h\in\Z}
        \Phi_i\pare{\frac hH}e\bigl(h(\alpha_i+(-1)^i\ell_i\alpha)\bigr)\right)d\alpha.
    \]
    Fourier inversion gives $\check\mu(bm)=a_{dmn''}$ for $m\sim N/(dn'')$. The estimates in the proof of~\cite[Theorem 1.7, equations (5.15)--(5.16)]{PascadiLargeSieve} give, for any $\eta>0$,
    \begin{align*}
        |\mu|(\R/\Z)&\ll_\eta H^\eta\pare{1+\frac{H}{L}},\\
        \mcI_{N/t}(\mu)&\ll_\eta H^{2\eta}L\min_iT_H(\alpha_i)\pare{1+\frac{H}{L}}^2.
    \end{align*}
    Translate $\mu$ by $t\xi/N$ to obtain $\mu_\xi$. Since $bt=dn''$, we have $\check\mu_\xi(bm)=a_{dmn''}e(dmn''\xi/N)$ for $m\sim N/(dn'')$, and $|\mu_\xi|(\R/\Z)=|\mu|(\R/\Z)$. Moreover,
    \begin{align*}
        T_{N/t}\pare{\alpha+\frac{t\xi}{N},\beta+\frac{t\xi}{N}}
        &\leq\pare{1+\frac{2N}{t}\left\|\frac{t\xi}{N}\right\|}T_{N/t}(\alpha,\beta)\\
        &\leq(1+2|\xi|)T_{N/t}(\alpha,\beta),
    \end{align*}
    so $\mcI_{N/t}(\mu_\xi)\leq(1+2|\xi|)\mcI_{N/t}(\mu)$.

    Apply Lemma~\ref{lemma: generalisation of pascadi large sieve} with the substitutions
    \[
        \begin{gathered}
            (N,a,\mu,\epsilon)\leftarrow\left(\frac{N}{dn''},d,\mu_\xi,\eta\right),\\
            a_m\leftarrow a_{dmn''}e\pare{\frac{dmn''\xi}{N}},\qquad
            F\equiv1,\qquad A\leftarrow\widetilde A:=C_\eta H^\eta A_{n''}.
        \end{gathered}
    \]
    For sufficiently large $C_\eta$, the required lower bound for $A$ holds: the coefficient norm is exactly $\norm{a_{nn''}}_2$, while $(d,q)\leq d$, $q\gg L^2$, and $N/(dn'')\leq N/t\ll HL$ give
    \begin{align*}
        \frac{(d,q)\bigl(N/(dn'')\bigr)^2}{q+N/n''}|\mu_\xi|(\R/\Z)^2
        &\ll_\eta H^{2\eta}\frac{N}{n''}\min\set{1,\frac HL}\pare{1+\frac{H}{L}}^2\\
        &\ll_\eta H^{2\eta}\frac{N}{n''}\pare{\frac{H}{L}+\frac{H^2}{L^2}}
        \ll_\eta H^{2\eta}A_{n''}^2.
    \end{align*}
    The condition on $X$ follows from
    \begin{align*}
        \frac{\widetilde A^2}{\mcI_{N/t}(\mu_\xi)}
        &\gg\frac{T_0NH}{n''L(H+L)\min_iT_H(\alpha_i)(1+|\xi|)},\\
        \max\set{1,\frac{\widetilde A^2}{\mcI_{N/t}(\mu_\xi)}}
        &\gg\frac{Y_{n''}}{1+|\xi|^2}.
    \end{align*}
    Since $d\cdot N/(dn'')=N/n''$ and $b\cdot N/(dn'')=N/t$, the large-sieve estimate applies in the required range
    \[
        X\ll\max\set{1,\frac{qn''}{N}}\frac{Y_{n''}}{1+|\xi|^2}.
    \]
    In this range $X\ll qN$. Taking $\eta$ sufficiently small in terms of the $\epsilon$ in Assumption~\ref{assumption: assumption on dispersion weights} absorbs the factor $H^{2\eta}$ and the large-sieve loss $\bigl((qN/n'')\max\{1,X\}\bigr)^{2\eta}$ into $(qNH)^\epsilon$. Together with~\eqref{eqn: dilated dispersion coefficients}, this proves the claim.
\end{proof}

For integers $m,n$ and $c\geq1$, the classical Kloosterman sum is
\[
    S(m,n;c):=\sum_{\substack{d\bmod c\\(d,c)=1}}e\pare{\frac{md+n\overline d}{c}}.
\]
The following proposition, which is a variant of~\cite[Corollary 5.11]{PascadiLargeSieve}, gives a bound for multilinear forms of Kloosterman sums with uniformity in the residue parameter $a$.
For sequences with two indices, we use the norm
\[
    \norm{b_{r,s}}_2^2:=\sum_{r,s}|b_{r,s}|^2.
\]
\begin{proposition}\label{propn: complete kloosterman bound}
    Let $R,S,M,N\geq1/2$, $C,Z\gg1$, $a\in\Z\setminus\{0\}$, $\epsilon>0$, and $\omega\in\R/\Z$. Suppose that:
    \begin{itemize}
        \item $(w_{r,s})$ is a complex sequence supported on $r\sim R$, $s\sim S$ with $(r,s)=(a,rs)=1$;
        \item for each such pair $(r,s)$, the function $\Phi_{r,s}:(0,\infty)^3\to\C$ is smooth, supported in $x,y,z\asymp1$, and satisfies
        \begin{align*}
            \abs{\frac{\partial^{j+k+\ell}}{\partial x^j\partial y^k\partial z^\ell}\Phi_{r, s}(x, y, z)}\ll_{j, k, \ell, \epsilon} Z^{(j+k)\epsilon}, \quad \forall j, k, \ell\geq 0;
        \end{align*}
        \item the tuple $(rs,n'',N,Z,(a_{n,r,s})_{n\sim N},A_{n'',r,s},Y_{n''})$ satisfies Assumption~\ref{assumption: assumption on dispersion weights} for each such pair and every integer $1\leq n''\leq2N$;
        \item $Y,\mcA>0$ are independent of $n''$ and, for every integer $1\leq n''\leq2N$, satisfy
        $$Y_{n''}\gg Y/n'',\quad\text{ and }\quad\norm{w_{r, s}A_{n'', r, s}}_2 \ll \mcA/\sqrt{n''}.$$
    \end{itemize}
    Then with a consistent choice of the $\pm$ sign, one has
    \begin{align*}
        \sum_{\substack{r\sim R\\s\sim S\\(r,s)=(a,rs)=1}}w_{r, s}\sum_{m\sim M}e(m\omega)\sum_{n\sim N}a_{n, r, s}\sum_{(c, r) = 1}\Phi_{r, s}\pare{\frac{m}{M},\frac{n}{N}, \frac{c}{C}}&S(m\overline{r}, \pm an; sc)\\ &\ll_{\epsilon} (\abs{a}RSMNCZ)^{O(\epsilon)}\mcA\cdot\mcL,
    \end{align*}
    where
    \begin{align*}
        \mcL:= \pare{\sqrt{\abs{a}}+\frac{C}{R\sqrt{SY}}}^{\theta}\sqrt{RSM}\pare{\frac{C^2}{R}(M+RS)(N+RS)+\abs{a}MN}^{1/2}.
    \end{align*}
\end{proposition}
\begin{proof}

     We closely follow the proof of \cite[Corollary 5.11]{PascadiLargeSieve}. Denote the sum on the left-hand side by $\mcS$. Let $Q:=RS$. By the divisor bound,
    \begin{align*}
        \mcS\ll Q^{o(1)}\sum_{\substack{Q<q\leq4Q\\(q,a)=1}}\max_{\substack{r\sim R\\ s\sim S\\ (r,s)=(a,rs)=1\\ rs = q}}\abs{w_{r, s}}\abs{\sum_{m\sim M}e(m\omega)\sum_{n\sim N}a_{n, r, s}\sum_{(c, r) = 1}\Phi_{r, s}\pare{\frac{m}{M}, \frac{n}{N}, \frac{c}{C}}S(m\overline{r}, \pm an; sc)},
    \end{align*}
    where we interpret any empty maximum as $0$. For each $q$ with a non-empty maximum, choose $r,s$ attaining the maximum. All subsequent sums run over these $q$. Put $A_{n'',q}=A_{n'',r,s}$, and let $w_q:=w_{r, s}$, $a_{n, q} = a_{n, r, s}$, $\Phi_{q}(x, y, z):= \Phi_{r, s}\pare{x, y, z(S/s)\sqrt{R/r}}$, and $\mathfrak{b}_q := 1/s$. Set $c'=s\sqrt r\,c$ and then rename $c'$ as $c$. This leaves us with bounding
    \begin{align*}
        Q^{o(1)}\sum_{\substack{Q<q\leq4Q\\(q,a)=1}}\abs{w_q}\abs{\sum_{m\sim M}e(m\omega)\sum_{n\sim N}a_{n, q}\sum_{c\in \mcC_{\infty \mathfrak{b}_q}}\Phi_{q}\pare{\frac{m}{M}, \frac{n}{N}, \frac{c}{S\sqrt{R}C}}S_{\infty\mathfrak{b}_q}(m, \pm an; c)},
    \end{align*}
    where the set $\mcC_{\infty \mathfrak{b}_q}$ and the generalised Kloosterman sum $S_{\infty\mathfrak{b}_q}(m, \pm an; c)$ are as defined in~\cite[Section 3.2]{PascadiLargeSieve}.

    Split $Q<q\leq4Q$ into the two dyadic intervals $(Q,2Q]$ and $(2Q,4Q]$. We treat either interval, relabelling its lower endpoint as $Q$, so that $Q\asymp RS$ throughout. By a standard Fourier transform argument, we see that
    \begin{align*}
        \mcS\ll_\epsilon Q^{o(1)}Z^{O(\epsilon)}S\sqrt{R}C\iint_{\R^2} \frac{\mcS(\zeta, \xi)d{\zeta}d{\xi}}{(1+\zeta^4)(1+\xi^4)},
    \end{align*}
    where
    \begin{align*}
        &\mcS(\zeta, \xi):=\\
        &\sum_{\substack{q\sim Q\\(q,a)=1}}\abs{w_q}\abs{\sum_{m\sim M}e\pare{m\pare{\omega+\frac{\zeta}{M}}}\sum_{n\sim N}a_{n, q}e\pare{n\frac{\xi}{N}}\sum_{c\in \mcC_{\infty\mathfrak{b}_q}}\frac{S_{\infty\mathfrak{b}_q}(m, \pm an; c)}{c}\varphi_{\zeta, \xi ,q}\pare{\frac{4\pi\sqrt{\abs{a}mn}}{c}}},
    \end{align*}
    and $\varphi_{\zeta, \xi, q}(z)$ are smooth functions supported in $z\asymp X^{-1}$, satisfying $\varphi^{(\ell)}_{\zeta, \xi, q}\ll_\ell X^{\ell}$ for
    \begin{align*}
        X:= \frac{S\sqrt{R}C}{\sqrt{\abs{a}MN}}.
    \end{align*}
    At this point, we apply the Kuznetsov trace formula. We only bound the contribution from the exceptional spectrum, as the contribution from the regular spectrum is easily bounded by applying the inverse Hecke relation followed by the ordinary spectral large sieve. Denote the contribution from the exceptional spectrum as $\mcS_{\mcM, \text{exc}}(\zeta, \xi)$. Indexing the chosen basis at level $q$ by $f_j$, we write $\lambda_j(q):=\lambda_{f_j}$, $\theta_j(q):=\theta_{f_j}$, $\kappa_j:=\kappa_{f_j}$, and $\rho_{j\mathfrak a}:=\rho_{f_j,\mathfrak a}$. We have
    \begin{align*}
        \mcS_{\mcM, \text{exc}}(\zeta, \xi)\ll
        \frac{1}{1+X^{-1}}\sum_{\substack{q\sim Q\\(q,a)=1}}\abs{w_q}\sum_{\lambda_j(q)<1/4}\frac{1+X^{\theta_j(q)}}{\cosh(\pi\kappa_j)}&\abs{\sum_{m\sim M}e\pare{m\pare{\omega+\frac{\zeta}{M}}}\overline{\rho_{j\infty_q}(m)}}\\&\times\abs{\sum_{n\sim N}a_{n, q}e\pare{n\frac{\xi}{N}}\rho_{j\mathfrak{b}_q}(an)}.
    \end{align*}
    Since $(a,q)=1$, applying Lemma~\ref{lemma: ABL factorisation} to the coefficients $a_{n,q}e(n\xi/N)$ gives
    \begin{align*}
        \mcS_{\mcM,\text{exc}}(\zeta,\xi)
        &\ll_\epsilon\abs{a}^\epsilon\sum_{\substack{n''\mid\abs{a}\\n''\leq2N}}
        \frac{(\abs{a}/n'')^{\theta/2}}{1+X^{-1}}\\
        &\quad\times\sum_{\substack{q\sim Q\\(q,a)=1}}\abs{w_q}
        \sum_{\lambda_j(q)<1/4}\frac{1+X^{\theta_j(q)}}{\cosh(\pi\kappa_j)}
        \abs{\sum_{n\sim N/n''}a_{nn'',q}e\pare{\frac{nn''}{N}\xi}\rho_{j\mathfrak{b}_q}(n)}\\
        &\quad\times\abs{\sum_{m\sim M}e\pare{m\pare{\omega+\frac{\zeta}{M}}}\overline{\rho_{j\infty}(m)}},
    \end{align*}
    Applying the Cauchy--Schwarz inequality, splitting $X= X_{0, n''}\sqrt{X_{1, n''}X_{2, n''}}$, we obtain
\begin{align*}
    \mcS_{\mcM, \text{exc}}(\zeta, \xi)\ll_\epsilon\abs{a}^\epsilon\sum_{\substack{n''\leq2N\\ n''\mid\abs{a}}}(\abs{a}/n'')^{\theta/2}\frac{(1+X_{0, n''})^{\theta}}{1+X^{-1}}\sqrt{\mcS_{M, n''}(\zeta)\mcS_{N, n''}(\xi)},
\end{align*}
where
\begin{align*}
    \mcS_{M, n''}(\zeta)&:= \sum_{\substack{q\sim Q\\(q,a)=1}}\sum_{\lambda_j(q)<1/4}\frac{1+X_{1, n''}^{\theta_j(q)}}{\cosh(\pi\kappa_j)}\abs{\sum_{m\sim M}e\pare{m\pare{\omega+\frac{\zeta}{M}}}\overline{\rho_{j\infty}(m)}}^2,\\
    \mcS_{N, n''}(\xi)&:=\sum_{\substack{q\sim Q\\(q,a)=1}}\abs{w_q}^2\sum_{\lambda_j(q)<1/4}\frac{1+X_{2, n''}^{\theta_j(q)}}{\cosh(\pi\kappa_j)}\abs{\sum_{n\sim N/n''}a_{nn'',q}e\pare{\frac{nn''}{N}\xi}\rho_{j\mathfrak{b}_q}(n)}^2.
\end{align*}
We set
$$X_{1, n''}:=\max\set{M, \frac{R^2S^2}{M}}$$
and bound $\mcS_{M, n''}(\zeta)$ as in the proof of~\cite[Corollary 5.11]{PascadiLargeSieve}. We set
\begin{align*}
    X_{2, n''}(\xi):= \max\set{1, \frac{Qn''}{N}}\frac{Y_{n''}}{1+\abs{\xi}^2}
\end{align*}
Since $q\asymp Q$, this choice is admissible in Assumption~\ref{assumption: assumption on dispersion weights}.

The selected pairs satisfy
\[
    \norm{w_qA_{n'',q}}_2^2
    :=\sum_{\substack{q\sim Q\\(q,a)=1}}|w_qA_{n'',q}|^2
    \leq\norm{w_{r,s}A_{n'',r,s}}_2^2\ll\mcA^2/n''.
\]
Substituting the bounds for $\mcS_{M, n''}$ and $\mcS_{N, n''}$, we obtain
\begin{align*}
    \mcS\ll_{\epsilon} (\abs{a}QMNCZ)^{O(\epsilon)}
    &\sum_{\substack{n''\mid\abs{a}\\n''\leq2N}}
    \left(\frac{\abs{a}}{n''}\right)^{\theta/2}
    \pare{1+\frac{X}{\sqrt{X_{1, n''}X_{2, n''}(0)}}}^{\theta}\frac{S\sqrt{R}C}{1+X^{-1}}\\
    &\times\pare{1+X^{-2}+\frac{M}{Q}}^{1/2}
    \pare{1+X^{-2}+\frac{N}{n''Q}}^{1/2}\sqrt{QM}\norm{w_qA_{n'',q}}_2.
\end{align*}
Since $X_{1,n''}\gg Q^2/M$ and $X_{2,n''}(0)\gg QY/N$, we have
\[
    \frac{X}{\sqrt{X_{1,n''}X_{2,n''}(0)}}
    \ll\frac{C}{R\sqrt{|a|SY}}.
\]
Also, $N/(n''Q)\leq N/Q$ and the norm hypothesis gives
\[
    \sum_{\substack{n''\mid|a|\\n''\leq2N}}
    (|a|/n'')^{\theta/2}\norm{w_qA_{n'',q}}_2
    \ll_\eta |a|^{\theta/2+\eta}\mcA.
\]
It remains to use the elementary inequality
\begin{align*}
    &\frac{S\sqrt R C}{1+X^{-1}}
    \left(1+X^{-2}+\frac M Q\right)^{1/2}
    \left(1+X^{-2}+\frac N Q\right)^{1/2}\\
    &\qquad\ll CS\sqrt R+\sqrt{|a|MN}
    +C\sqrt{SM}+C\sqrt{SN}+C\sqrt{MN/R}\\
    &\qquad\ll\left(\frac{C^2}{R}(M+RS)(N+RS)+|a|MN\right)^{1/2}.
\end{align*}
Here the first bound follows by expanding the product
\[
    \frac{(CS\sqrt R+\sqrt{|a|MN}+C\sqrt{SM})
    (CS\sqrt R+\sqrt{|a|MN}+C\sqrt{SN})}
    {CS\sqrt R+\sqrt{|a|MN}}.
\]
Combining these estimates, and taking $\eta$ sufficiently small in terms of $\epsilon$, proves the proposition.
\end{proof}
The following proposition will serve as a substitute for~\cite[Proposition 3.7]{PascadiLvlDist}. It incorporates the uniformity in the residue parameter $a$ that we require for our application to bound multilinear forms of incomplete Kloosterman sums.
\begin{proposition}\label{propn: incomplete kloosterman sums}
    Let $R,S,N\geq1/2$, $C,D,Z\gg1$, $a\in\Z\setminus\{0\}$, and $\epsilon>0$. Suppose that:
    \begin{itemize}
        \item $(w_{r,s})$ is a complex sequence supported on $r\sim R$, $s\sim S$ with $(r,s)=(a,rs)=1$;
        \item for each such pair $(r,s)$, the function $\Phi_{r,s}:(0,\infty)^3\to\C$ is smooth, supported in $x,y,z\asymp1$, and satisfies
        \begin{align*}
            \abs{\frac{\partial^{j+k+\ell}}{\partial x^j\partial y^k\partial z^\ell}\Phi_{r, s}(x, y, z)}\ll_{j, k, \ell, \epsilon} Z^{j\epsilon}, \quad \forall j, k, \ell\geq 0;
        \end{align*}
        \item the tuple $(rs,n'',N,Z,(a_{n,r,s})_{n\sim N},A_{n'',r,s},Y_{n''})$ satisfies Assumption~\ref{assumption: assumption on dispersion weights} for each such pair and every integer $1\leq n''\leq2N$;
        \item $Y,\mcA>0$ are independent of $n''$ and, for every integer $1\leq n''\leq2N$, satisfy
        $$Y_{n''}\gg Y/n'',\quad\text{ and }\quad\norm{w_{r, s}A_{n'', r, s}}_2 \ll \mcA/\sqrt{n''}.$$
    \end{itemize}
    Then with a consistent choice of the sign $\pm$, we have
    \begin{align*}
        \sum_{\substack{r\sim R\\ s\sim S\\ (r,s)=(a,rs)=1}}w_{r, s}\sum_{n\sim N}a_{n, r, s}\sum_{\substack{c, d\\ (rd, sc)=1}}e\pare{\pm an\frac{\overline{rd}}{sc}}\Phi_{r, s}\pare{\frac{n}{N}, \frac{c}{C}, \frac{d}{D}}\\
        \ll_{\epsilon} (\abs{a}RSNCDZ)^{O(\epsilon)}\mcA\cdot \mathcal{I},
    \end{align*}
    where
    \begin{align*}
        \mcI^2 := D^2NR + \pare{\abs{a}+\frac{C^2}{R^2SY}}^{\theta}(CS(C+DR)(RS+N) +\abs{a}RSN).
    \end{align*}
\end{proposition}
\begin{proof}
    This is a standard deduction from Proposition~\ref{propn: complete kloosterman bound}. See the proof of~\cite[Corollary 5.14]{PascadiLargeSieve} for the details.
\end{proof}
The following lemma gives a bound for multilinear forms of incomplete Kloosterman sums with coefficients specialised to our application. It plays a similar role as~\cite[Lemma 6.1]{PascadiLvlDist}, but with the uniformity in the residue parameter $a$ that we require for our application.
\begin{lemma}\label{lemma: multilinear kloosterman}
    Let $x\geq2$, $\epsilon>0$, $1\ll N,T,H,K,L\ll x$ with $TH\ll N$, and $a\in\Z\setminus\{0\}$, $b,d\in\Z_{>0}$ with $1\leq\abs{a}\leq x^{1+\epsilon}$, $1\leq d\leq x^{2\epsilon}$, and $(a,d)=1$. Let $\Phi_i(t)$ be smooth functions supported in $t\asymp1$ with $\Phi_i^{(j)}\ll_j1$, and set
    \begin{align*}
        \phi(h_1, h_2):=\Phi_1\pare{\frac{h_1}{H}}\Phi_2\pare{\frac{h_2}{H}}e(h_1\alpha_1+h_2\alpha_2),
    \end{align*}
    where $\alpha_i\in \R/\Z$ have $\min_i T_H(\alpha_i)\ll x^{2\epsilon}$. Then for any smooth function $\Phi:(0,\infty)^2\to\C$ supported in $x_i\asymp1$, satisfying $\partial_{x_1}^{j_1}\partial_{x_2}^{j_2}\Phi(x_1,x_2)\ll_{j_1,j_2,\epsilon}1$ for all $j_1,j_2\geq0$, one has, uniformly in $b$,
    \begin{align*}
        &\sum_{\substack{n,n'\sim N\\(nn',a)=1}}\abs{\sum_{\substack{1\leq \abs{t}\leq T\\(bt,nn')=1\\bt\mid n-n'}}\sum_{\substack{h, h'\\ e = t(n'h-nh')\neq 0}}\phi(h, h')\sum_{\substack{k, \ell\\(k, dnn'\ell)=1}}\Phi\pare{\frac{\ell}{L}, \frac{k}{K}}e\pare{\frac{ae\overline{dnn'\ell}}{k}}}\\
        &\ll_{\epsilon} x^{6\epsilon}THN\pare{L^2THN^3+\pare{\abs{a}+\frac{K^2}{N^3TH^2}}^{\theta}(KN^4L+\abs{a}N^3TH+K^2N^2)}^{1/2}.
    \end{align*}
\end{lemma}
\begin{proof}
    We follow the proof of~\cite[Lemma 6.1]{PascadiLvlDist}, applying Proposition~\ref{propn: incomplete kloosterman sums} in place of~\cite[Proposition 3.7]{PascadiLvlDist} to retain uniformity in $a$. Denote the sum in question by $\mcK$ and split it into $\mcK(n=n')$ and $\mcK(n\neq n')$.

    For the diagonal, set $q=dn$. For each sign and dyadic range $e\sim E\ll TH$, put $Q=dN$ and, for $q\sim Q$ with $(q,a)=1$, define
    \[
        a_{e,q}:=\mathbf{1}_{d\mid q}
        \sum_{\substack{1\leq|t|\leq T\\(bt,q/d)=1}}
        \sum_{\substack{h,h'\\\pm t(h-h')=e}}\phi(h,h').
    \]
    Thus $a_{e,q}=0$ when $d\nmid q$. With
    \[
        \mcK_1(E,Q):=\sum_{\substack{q\sim Q\\(q,a)=1}}
        \abs{\sum_{e\sim E}a_{e,q}
        \sum_{\substack{k,\ell\\(k,q\ell)=1}}
        \Phi\pare{\frac{\ell}{L},\frac{k}{K}}
        e\pare{\frac{\pm ae\overline{q\ell}}{k}}},
    \]
    we have $\mcK(n=n')\ll x^{o(1)}\sup_{E\ll TH}\mcK_1(E,Q)$, with both signs treated in the same way. Since $t\mid e$, the divisor bound gives $|a_{e,q}|\ll x^{o(1)}H$, and hence
    \[
        \sum_{\substack{q\sim Q\\(q,a)=1}}\sum_{e\sim E/e''}|a_{e''e,q}|^2
        \ll x^{o(1)}\frac{H^2EQ}{e''}
        \qquad(1\leq e''\leq2E).
    \]
    For each $q$ and $1\leq e''\leq2E$, Lemma~\ref{lemma: assumption for dispersion weights}(i) shows that $(q,e'',E,1,(a_{e,q})_{e\sim E},A_{e'',q},1)$ satisfies Assumption~\ref{assumption: assumption on dispersion weights}, where
    \[
        A_{e'',q}:=\left(\sum_{e\sim E/e''}|a_{e''e,q}|^2\right)^{1/2}.
    \]
    In particular, we may take $Z=1$ and $Y_{e''}=Y=1$. The preceding norm bound gives the remaining hypothesis of Proposition~\ref{propn: incomplete kloosterman sums}:
    \[
        \left(\sum_{\substack{q\sim Q\\(q,a)=1}}A_{e'',q}^2\right)^{1/2}
        \ll\frac{\mcA}{\sqrt{e''}},\qquad
        \mcA=x^{o(1)}H\sqrt{EQ}.
    \]
    Applying Proposition~\ref{propn: incomplete kloosterman sums} with $R=Q$, $S=1/2$, $N=E$, $C=K$, and $D=L$, we obtain
    \begin{align*}
        \mcK_1(E,Q)\ll_\epsilon x^{\epsilon/2}H\sqrt{EQ}\Bigg(&L^2EQ
        +\pare{|a|+\frac{K^2}{Q^2}}^\theta\\
        &\times\bigl(K(K+LQ)(Q+E)+|a|EQ\bigr)\Bigg)^{1/2}.
    \end{align*}
    This bound is non-decreasing in $E$ and $Q$. Using $E\ll TH$, $Q=dN$, and $d\leq x^{2\epsilon}$ therefore gives
    \begin{align*}
        \mcK(n=n')\ll_\epsilon x^{4\epsilon}H\sqrt{THN}\Bigg(&L^2THN
        +\pare{|a|+\frac{K^2}{N^2}}^\theta\\
        &\times\bigl(K(K+LN)(N+TH)+|a|NTH\bigr)\Bigg)^{1/2}.
    \end{align*}
    The above diagonal contribution $\mcK(n=n')$ is acceptable, as $TH\ll N$ and $\theta\leq1/3$ give
    \begin{align*}
        \frac{H}{TN}\pare{|a|+\frac{K^2}{N^2}}^\theta
        &\ll\frac{H}{TN}(NTH^2)^\theta
        \pare{|a|+\frac{K^2}{N^3TH^2}}^\theta\\
        &\ll N\pare{|a|+\frac{K^2}{N^3TH^2}}^\theta.
    \end{align*}
    Furthermore, we have
    \[
        N\bigl(K(K+LN)(N+TH)+|a|NTH\bigr)
        \ll K^2N^2+KLN^4+|a|N^3TH.
    \]

    For the off-diagonal terms, let $n_0:=(n,n')$, replace $(n,n')$ by $(n_0n,n_0n')$, and put $e:=t(n'h-nh')$ in the new variables. Applying the triangle inequality in $t$ and then dropping the additional restrictions involving $b$ gives
    \begin{align*}
        \mcK(n\neq n')\ll \sum_{\substack{n_0\leq 2N\\ n, n'\sim N/n_0\\(n,n')=1\\(a,n_0nn')=1}}\sum_{\substack{1\leq\abs{t}\leq T\\t\mid n-n'\neq0\\(t,n_0nn')=1}}\abs{\sum_{\substack{h, h'\\ e = t(n'h-nh')\neq 0}}\phi(h, h')\sum_{\substack{k, \ell\\ (k,dn_0nn'\ell)=1}}\Phi\pare{\frac{\ell}{L}, \frac{k}{K}}e\pare{\frac{ae\overline{dn_0nn'\ell}}{k}}}
    \end{align*}
    The support conditions in the expression above, together with $(a,d)=1$, give $(a,dn_0nn')=1$. Splitting by the sign of $e$, we put $n_0$, $\abs{e}$, and $q=dn_0nn'$ in dyadic ranges and use the divisor bound to write
    \begin{align*}
        \mcK(n\neq n') \ll x^{o(1)}\sup_{\substack{N_0\ll N\\ E\ll THN/N_0\\ Q\asymp dN^2/N_0}}\mcK_2(N_0, E, Q),
    \end{align*}
    where
    \begin{align*}
        \mcK_2&:= \sum_{\substack{q\sim Q\\(q,a)=1}}\max_{\substack{n_0\sim N_0,\ n,n'\sim N/n_0\\q=dn_0nn',\ (n,n')=1\\1\leq\abs{t}\leq T,\ t\mid n-n'\neq0\\(t,n_0nn')=1}}\abs{\sum_{e\sim E}\sum_{\substack{h, h'\\ t(n'h-nh')=\pm e}}\phi(h, h')\sum_{\substack{k, \ell\\ (k, q\ell) = 1}}\Phi\pare{\frac{\ell}{L},\frac{k}{K}}e\pare{\frac{\pm ae\overline{q\ell}}{k}}}.
    \end{align*}
    For a nonempty maximum, choose a maximising tuple and set
    \begin{align*}
        a_{e, q}:= \sum_{\substack{h, h'\in \Z\\ \pm t(q)(n'(q)h-n(q)h')=e}}\phi(h, h'),
    \end{align*}
    with $a_{e,q}=0$ for an empty maximum. We obtain
    \begin{align*}
        \mcK_2 = \sum_{\substack{q\sim Q\\(q,a)=1}}\abs{\sum_{e\sim E}a_{e,q}\sum_{\substack{k, \ell\\ (k, q\ell) = 1}}\Phi\pare{\frac{\ell}{L},\frac{k}{K}}e\pare{\frac{\pm ae\overline{q\ell}}{k}}}.
    \end{align*}
    Apply Lemma~\ref{lemma: assumption for dispersion weights}(ii) with
    \[
        N\leftarrow E,\quad n''\leftarrow e'',\quad L\leftarrow N/N_0,\quad
        (\ell_1,\ell_2)\leftarrow(n'(q),n(q)),\quad t\leftarrow\pm t(q),\quad T_0\leftarrow T.
    \]
    Its hypotheses hold since $(n(q),n'(q))=1$ and $q=dn_0n(q)n'(q)\gg(N/N_0)^2$. Thus, for every integer $1\leq e''\leq2E$, the tuple $(q,e'',E,x,(a_{e,q})_{e\sim E},A_{e'',q},Y_{e''})$ satisfies Assumption~\ref{assumption: assumption on dispersion weights}, where
    \begin{align*}
        Y&:=\frac{EH}{(H+N/N_0)(N/N_0)\min_iT_H(\alpha_i)}
        \gg_\epsilon x^{-2\epsilon}\frac{EHN_0}{(H+N)N},\\
        Y_{e''}&:=\max\set{1,\frac{Y}{e''}}\geq\frac{Y}{e''},\\
        A_{e'',q}&:=\pare{\sum_{e\sim E/e''}\abs{a_{e''e,q}}^2}^{1/2}
        +\sqrt{\frac{TE}{e''}}\sqrt{\frac{HN_0}{N}+\frac{H^2N_0^2}{N^2}}.
    \end{align*}
    Up to an absolute constant, $A_{e'',q}$ plays the role of $A_{n''}$ in Lemma~\ref{lemma: assumption for dispersion weights}(ii). Since $H\ll x$, we may replace $Z=H$ by $Z=x$ in Assumption~\ref{assumption: assumption on dispersion weights}.
    Since $(n_1,n_2)=1$, any two solutions of $h_1n_1-h_2n_2=\pm e/t$ have $h_1$ congruent modulo $n_2$. The divisor bound therefore gives
    \begin{align*}
        \sum_{\substack{q\sim Q\\(q,a)=1}}\sum_{\substack{e\sim E\\e''\mid e}}\abs{a_{e,q}}^2
        &\ll\sum_{\substack{e\sim E\\e''\mid e}}\sum_{\substack{1\leq\abs{t}\leq T\\t\mid e}}
        \sum_{n_0\sim N_0}\sum_{\substack{n_1,n_2\sim N/n_0\\(n_1,n_2)=1}}
        \left(\sum_{\substack{h_1,h_2\asymp H\\h_1n_1-h_2n_2=\pm e/t}}1\right)^2\\
        &\ll_\epsilon x^{2\epsilon}\frac{TE}{e''}NH\left(1+\frac{HN_0}{N}\right)
        =x^{2\epsilon}\frac{TE}{e''}(HN+H^2N_0).
    \end{align*}
    Additionally,
    \begin{align*}
        \sum_{\substack{q\sim Q\\(q,a)=1}}\frac{TE}{e''}\pare{\frac{HN_0}{N}+\frac{H^2N_0^2}{N^2}}\ll x^{2\epsilon}\frac{TE}{e''}(HN+H^2N_0).
    \end{align*}
    Thus,
    \begin{align*}
        \norm{A_{e'', q}}_2^2 \ll x^{2\epsilon}\frac{TE}{e''}(HN+H^2N_0).
    \end{align*}
    Since every contributing level satisfies $(a,q)=1$, Proposition~\ref{propn: incomplete kloosterman sums} gives
    \begin{align*}
        \mcK_2\ll_{\epsilon} x^{4\epsilon}\sqrt{TE(HN+H^2N_0)}\Bigg(L^2EQ+&\pare{{\abs{a}}+\frac{K^2}{Q^2 \frac{EHN_0}{(H+N)N}}}^{\theta }\\&(K(K+LQ)(Q+E) +\abs{a}EQ)\Bigg)^{1/2}.
    \end{align*}
    After using $E\ll THN/N_0$ and $Q\asymp dN^2/N_0$, the resulting upper bound is non-increasing in $N_0\geq1$, since $\theta\leq1/2$. We may therefore take $N_0=1$; the factors involving $d\leq x^{2\epsilon}$ are absorbed into the stated $x^{6\epsilon}$ factor.
\end{proof}

We now obtain a power-saving bound for an exponential sum, similar to the estimates in~\cite[Lemma~6.2]{PascadiLvlDist} and~\cite[Proposition~1]{Drappeau}. We use the smooth-weight formulation of Pascadi and keep the dependence on the residue parameter $a$ explicit.
\begin{lemma}\label{lemma: exponential sum bound for convolution}
    Let $x\geq2$, $\epsilon>0$ be small enough, $a\in \Z\setminus \set{0}$, $v,d_1,d_2\in\Z_{>0}$ with $\abs{a}\ll x^{1+\epsilon}$ and $(a,vd_1d_2)=1$, $\theta:=7/32$, and $1\ll M, K, N, L, H, R\ll x$ satisfy
    \begin{equation}\label{eqn: conditions on parameters to bound exponential sum}
        \begin{aligned}
            &\abs{vd_1d_2}\ll x^\epsilon, \quad NL\ll x^\epsilon K, \quad R\ll K\ll \min\set{x^{-3\epsilon}MN, LN^2},\\
            &H\ll x^{\epsilon}\frac{R}{M},\quad K\ll x^{-25\epsilon}\sqrt{MNR}, \quad K^{3+\theta}N^{2-3\theta}\ll x^{-200\epsilon}M^{2-2\theta}R^{2+\theta}L, \\
            & K^{3}N^2\abs{a}^{\theta}\ll x^{-100\epsilon} M^2R^{2}L, \quad \abs{a}\ll x^{\epsilon}MLN.
        \end{aligned}
    \end{equation}
    Let $(u_k)_{K<k\leq 4K}, (\beta_n)_{n\sim N}, (\lambda_{\ell})_{\ell \sim L}$ be complex sequences such that $\abs{u_k}\leq \tau(k)$, $\abs{\beta_n}\leq 1$, $\abs{\lambda_{\ell}}\leq 1$, and
    \begin{align*}
        (k,avd_1d_2)>1\implies u_k = 0, \quad (n\ell,avd_1)>1\implies \beta_n\lambda_{\ell} = 0.
    \end{align*}
    Then for any smooth functions $\Phi(t), \Psi(t)$ supported in $t\asymp 1$ with $\Phi^{(j)}, \Psi^{(j)}\ll_j 1$, and any $\omega\in \R/\Z$ with $T_H(\omega)\ll x^\epsilon$, one has
    \begin{align*}
        \sum_{\substack{r\sim R\\ (r, avd_1d_2) = 1}}\frac{M}{r}\Psi\pare{\frac{r}{R}}\sum_{\substack{k, n, \ell\\ d_1k\equiv d_2n\ell \pmod r\\ (d_1k, d_2n\ell) = 1}}u_k\beta_n\lambda_\ell \sum_{h\in \Z}e(h\omega)\Phi\pare{\frac{h}{H}}e\pare{\frac{-ha\overline{vd_1d_2k}}{r}}\ll_\epsilon \frac{KMNL}{x^{10\epsilon}R}.
    \end{align*}
\end{lemma}
\begin{proof}
    We follow the proof of~\cite[Proposition 1]{Drappeau} and the smooth $h$-sum reduction in~\cite[Lemma 6.2]{PascadiLvlDist}. Denote the sum by $\mcR$, and put
    \[
        \nu:=vd_1d_2,\qquad T:=4\frac{\max\{d_1K,d_2NL\}}{R}
        \ll x^{2\epsilon}\frac{K}{R}.
    \]
    In the reciprocity identity~\cite[(3.20)]{Drappeau}, we retain the factor $e(-ah/(\nu kr))$: the error from discarding it in~\cite[(3.21)]{Drappeau} is not admissible throughout our range of $a$. Set
    \[
        D:=1+\frac{|a|H}{\nu KR},\qquad
        \nu D=\nu+\frac{|a|H}{KR}
        \ll x^{3\epsilon}.
    \]

    For fixed $k,\ell,t$, write $r(n)=(d_2n\ell-d_1k)/t$ and $b(n)=aH/(\nu kr(n))$. With a fixed nonnegative smooth cutoff $\Phi_0$ supported in $u\asymp1$ and equal to $1$ on $[1,4]$ and on the support of $\Phi$, Fourier inversion gives
    \[
        \Phi\pare{\frac hH}e\pare{-\frac{ah}{\nu kr(n)}}
        =\Phi_0\pare{\frac hH}\int_\R\widehat\Phi(s+b(n))e\pare{\frac{sh}{H}}\,ds.
    \]
    Let $[u_0,u_1]$ be the interval defined by $N\leq u\leq2N$ and $R\leq r(u)\leq2R$, and put $W_s(u):=(R/r(u))\Psi(r(u)/R)\widehat\Phi(s+b(u))$. Since $b(u)$ is monotone and $|b(u)|\ll D$, differentiation and the rapid decay of $\widehat\Phi$ give
    \[
        \int_\R\sup_{k,\ell,t}\left(|W_s(u_1)|+\int_{u_0}^{u_1}|W_s'(u)|\,du\right)ds
        \ll D,
    \]
    where the supremum is over the contributing triples $(k,\ell,t)$.
    For $\rho\in\R/\Z$ with $\nu\rho=0$, put
    \[
        d:=\nu d_2,\qquad \alpha_1=-\alpha_2=\omega+\rho+s/H.
    \]
    Partial summation in $n$, followed by the interval separation, character decomposition and Cauchy--Schwarz argument in~\cite[Section 3.5, (3.22)]{Drappeau}, now gives the analogue of~\cite[(6.7)]{PascadiLvlDist}:
    \[
        \mcR\ll_\epsilon x^{4\epsilon}\frac{M}{R}\sqrt{KLT}
        \max_{\sigma\mid |a|}
        \sup_{\substack{|s|\leq x^\epsilon D\\\rho\in\R/\Z,\ \nu\rho=0}}
        \mcB_\sigma(s,\rho)^{1/2}+O(x^{-100}),
    \]
    where
    \[
        \mcB_\sigma(s,\rho):=\sum_{\substack{n,n'\sim N\\(nn',a)=1}}
        \left|\begin{aligned}
            &\sum_{\substack{1\leq|t|\leq T\\(nn',vd_1\sigma t)=1\\vd_1\sigma t\mid n-n'}}
            \sum_{h,h'\in\Z}
            \Phi_0\pare{\frac hH}\Phi_0\pare{\frac{h'}H}
            e(h\alpha_1+h'\alpha_2)\\
            &\quad\times\sum_{\substack{k,\ell\\(k,dnn'\ell)=1}}
            \Phi_0\pare{\frac kK}\Phi_0\pare{\frac\ell L}
            e\pare{\frac{at(n'h-nh')\overline{dnn'\ell}}{k}}
        \end{aligned}\right|.
    \]
    The divisor $\sigma\mid |a|$ arises from detecting $(r,a)=1$. The tails in $s$ are negligible by the rapid decay of $\widehat\Phi$, and the divisor and interval-separation losses are absorbed into $x^\epsilon\nu D\ll x^{4\epsilon}$.
    The support conditions and $(a,\nu)=1$ give $(a,dnn')=1$. Also, $TH\ll N$ and
    \[
        T_H(\alpha_1)
        \leq\nu(1+2|s|)T_H(\omega)
        \ll x^{2\epsilon}\nu D\ll x^{5\epsilon}.
    \]
    The terms $n'h=nh'$ contribute $\ll_\epsilon x^\epsilon KLNHT$ by~\cite[(3.24)]{Drappeau}. For the remaining terms, apply Lemma~\ref{lemma: multilinear kloosterman} with $b=vd_1\sigma$ and $5\epsilon/2$ in place of $\epsilon$. Using $|a|N^3TH\ll x^{4\epsilon}KLN^4$, and then simplifying as in the proof of~\cite[Lemma 6.2]{PascadiLvlDist}, gives
    \[
        \mcR\ll_\epsilon x^{7\epsilon}\frac{K^2L\sqrt{MN}}{R^{3/2}}
        +x^{15\epsilon}\frac{K^{7/4}L^{3/4}M^{1/2}N^{3/2}}{R^{3/2}}
        \pare{|a|+\frac{KM^2}{N^3R}}^{\theta/4}+O(x^{-100}).
    \]
    The hypotheses~\eqref{eqn: conditions on parameters to bound exponential sum} bound this by $x^{-10\epsilon}KMNL/R$, as required.
\end{proof}
We now give the proof of Proposition~\ref{propn: triple convolution with a uniformity}.
\begin{proof}[Proof of Proposition~\ref{propn: triple convolution with a uniformity}]
    We first prove the result in the special case where the coefficients $\beta_n$ and $\gamma_\ell$ satisfy
    \begin{equation}\label{eqn: coprime support for triple convolution}
        (n\ell,a_1)>1\implies\beta_n\gamma_\ell=0.
    \end{equation}
    We will remove this auxiliary assumption after proving the result in this case.
    We closely follow the proofs of~\cite[Théorème 3]{Drappeau} and~\cite[Proposition 6.3]{PascadiLvlDist}. Setting $u_k:=\sum_{n\ell=k}\beta_n\gamma_\ell$ (so $u_k=0$ when $(k,a_1)>1$), $K:=NL$, and following the steps in~\cite[Section 3.1]{Drappeau} with $\omega_\epsilon$ replaced by $\omega_D$, one is reduced to estimating the sums $\mcS_1, \mcS_2,$ and $\mcS_3$ defined as follows:
    \[
\begin{aligned}
&\mathcal{S}_1 :=
\sum_{\substack{R < r \leq 2R \\ (r, a_1 a_2) = 1}}
\ \sum_{\substack{(m, r) = 1}} f(m)
\left|
    \sum_{\substack{k \equiv a_1  \overline{a_2m} \pmod r}} u_k
\right|^2, \\[3ex]
&\mathcal{S}_2 :=
\sum_{\substack{R < r \leq 2R \\ (r, a_1 a_2) = 1}}
\frac{1}{\varphi(r)}
\sum_{\substack{(m, r) = 1}} f(m)
\mathop{\sum\sum}_{\substack{(k_1, r) = 1\\ k_2 \equiv a_1  \overline{a_2m} \pmod r}}
u_{k_1} \omega_{D}(mk_1 \overline{a_1 }a_2; r) \overline{u_{k_2}}, \\[3ex]
&\mathcal{S}_3 :=
\sum_{\substack{R < r \leq 2R \\ (r, a_1 a_2) = 1}}
\frac{1}{\varphi(r)^2}
\sum_{\substack{(m, r) = 1}} f(m)
\left|
    \sum_{\substack{(k, r) = 1}} u_k \omega_{D}(mk \overline{a_1 }a_2; r)
\right|^2,
\end{aligned}
\]
where $f(m):=\Phi_0(m/M)$ for a $C^{\infty}$-function $\Phi_0:\R_{>0}\to \R_{\geq0}$ which is compactly supported on $[1/2, 3]$ and is identically $1$ in the interval $[1, 2]$. Let $\mcD$ denote the left-hand side of the proposition. Since $|\alpha_m|\leq1$, Cauchy--Schwarz gives
\begin{equation}\label{eqn: dispersion cauchy schwarz}
    \mcD^2\ll MR(\mcS_1-2\Re\mcS_2+\mcS_3).
\end{equation}
Thus it suffices to estimate this combination of the three sums.

    The sum $\mcS_3$ can be estimated as in ~\cite[Section 3.2]{Drappeau} as the argument provided there is uniform in the residue class $a_1$, and we have $D\leq x^{\epsilon}$. Thus, we have
    \begin{align*}
        \mcS_3 = \widehat{f}(0)X_3+O_{\epsilon}\pare{x^{1-\delta}\frac{K}{R}},
    \end{align*}
    where
    \begin{align*}
        X_3:= \sum_{\substack{r\sim R\\(r, a_1a_2) = 1}}\frac{1}{r}\sum_{\substack{0<b<r\\ (b, r)=1}}\abs{\frac{1}{\phi(r)}\sum_{(k, r) = 1}u_k\omega_{D}(k\overline{b}; r)}^2.
    \end{align*}

    We note that the sum $\mcS_2$ can be estimated as in ~\cite[Section 3.3]{Drappeau} as the argument goes through even when $\abs{a_1}\leq x^{1+\epsilon}$ instead of $\abs{a_1}\leq x$, and we have $D\leq x^{\epsilon}$. Thus, we have
    \begin{align*}
        \mcS_2 = \widehat{f}(0)X_2+O_{\epsilon}\pare{x^{1-\delta}\frac{K}{R}},
    \end{align*}
    where
     \begin{align*}
        X_2:= \sum_{\substack{r\sim R\\(r, a_1a_2) = 1}}\frac{1}{r\phi(r)}\mathop{\sum\sum}_{\substack{(k_1, r) = 1\\(k_2, r) = 1}}u_{k_1}\overline{u_{k_2}}\omega_{D}(k_1\overline{k_2}; r).
     \end{align*}

    We write $d\mid n^\infty$ if every prime divisor of $d$ divides $n$. We retain coprimality to $a_1$ in Drappeau's decomposition~\cite[Section 3.4]{Drappeau}. Indeed, in its notation,
    \[
        v=(k_1,k_2),\qquad d_1\mid v^\infty,\qquad
        e_1\mid a_2^\infty,\qquad e_2\mid a_2,\qquad
        d_2=(k_2/(ve_2),v),
    \]
    so $(a_1,vd_1e_1e_2d_2)=1$. Dividing the coefficient indices by these factors preserves their coprimality to $a_1$. Under the substitution in that decomposition,
    \[
        (v,d_1,d_2)\longleftarrow
        \left(\frac{\abs{a_2}v}{e_2d_2},e_1d_1,e_2d_2\right),
        \qquad (vd_1d_2)_{\mathrm{new}}=\abs{a_2}vd_1e_1,
    \]
    the new parameters remain coprime to $a_1$; the sign of $a_2$ is absorbed into the residue parameter. Thus the additional coprimality and support hypotheses of Lemma~\ref{lemma: exponential sum bound for convolution} hold.

    To bound $\mcS_1$, we follow the decomposition in~\cite[proof of Proposition 6.3]{PascadiLvlDist}, applying Lemma~\ref{lemma: exponential sum bound for convolution} (with $\eta>0$ in place of $\epsilon>0$) in place of~\cite[Lemma 6.2]{PascadiLvlDist}. The preliminary steps are uniform in the residue class $a_1$, and $D\leq x^\epsilon$. Choose $0<\delta\ll\kappa\ll\eta\ll\epsilon$, using $x^\kappa$ for the decomposition cutoffs. This allows us to verify that the conditions in~(\ref{eqn: hypothesis for triple convolution}) and $\abs{a_1}\leq x^{1+\delta}$ imply the conditions of~(\ref{eqn: conditions on parameters to bound exponential sum}). Thus, we have
    \begin{align*}
        \mcS_1 = \widehat{f}(0)X_1+O_{\epsilon}\pare{x^{1-\delta}\frac{K}{R}},
    \end{align*}
    where
    \begin{align*}
        X_1:= \sum_{\substack{r\sim R\\(r, a_1a_2) = 1}}\frac{1}{r}\sum_{\substack{0<b<r\\ (b, r)=1}}\abs{\sum_{k\equiv b\pmod r}u_k}^2.
    \end{align*}

It remains to bound the contribution from the main terms. For this, we follow the steps in~\cite[Section 3.6]{Drappeau}. Doing so, we obtain
    \[
\begin{aligned}
X_1 - 2\Re (X_2) + X_3 &= \sum_{\substack{R < r \leq 2R \\ (r,a_1 a_2)=1}} \frac{1}{r} \sum_{\substack{0 < b < r \\ (b,r) = 1}}  \left| \sum_{k \equiv b \, \pmod r} u_k - \frac{1}{\varphi(r)} \sum_{(k, r) = 1} u_k \omega_D(k\overline{b}; r) \right|^2 \\
&\leq \frac{1}{R} \sum_{\substack{R < r \leq 2R \\ (r,a_1 a_2)=1}} \frac{1}{\varphi(r)} \sum_{\substack{\chi \text{ prim.} \\ \mathrm{cond}(\chi)> D \\ \mathrm{cond}(\chi) \mid r}} \left| \sum_{(k, r) = 1} u_k \chi(k) \right|^2.
\end{aligned}
\]
Rearranging the character sums, using Cauchy--Schwarz, and applying partial summation, we obtain
\begin{align*}
X_1 - 2\Re (X_2) + X_3 &\leq \frac{1}{R} \sum_{D < s \leq 2R} \sum_{\substack{\chi \pmod s \\ \chi \text{ prim.}}} \sum_{\substack{r \leq 2R \\ s \mid r}}\frac{1}{\varphi(r)} \left| \sum_{(k, r/s)=1} u_{k} \chi(k) \right|^2 \displaybreak[1]\\
&\leq \frac{1}{R} \sum_{D < s \leq 2R} \sum_{\substack{\chi \pmod s \\ \chi \text{ prim.}}}
 \sum_{\substack{r \leq 2R \\ s \mid r}} \frac{\tau(r/s)}{\varphi(r)} \sum_{d \mid r/s} \left| \sum_{d \mid k} u_{k} \chi(k) \right|^2 \displaybreak[1]\\
 &\ll \frac{(\log R)^2}{R} \sum_{d \leq R} \frac{\tau(d)}{\varphi(d)} \sum_{D < s \leq 2R} \frac{1}{\varphi(s)} \sum_{\substack{\chi \pmod s \\ \chi\text{ prim.}}} \left| \sum_{k'} u_{d k'} \chi(k') \right|^2\\
 &= \frac{(\log R)^2}{R} \sum_{d \leq R} \frac{\tau(d)}{\varphi(d)} \int_{D}^\infty \sum_{D < s \leq \min\{2R, t\}} \frac{s}{\varphi(s)} \sum_{\substack{\chi \pmod s \\ \chi\text{ prim.}}} \left| \sum_{k'} u_{d k'} \chi(k') \right|^2 \frac{dt}{t^2}
\end{align*}
Here $s>D\geq1$ is an integer, so $s\geq2$; hence $d\mid r/s$ and $r\leq2R$ imply $d\leq R$. Applying the multiplicative large sieve~\cite[Lemme 6]{Drappeau} and using $\abs{u_k}\leq\tau(k)$, we get
\begin{align*}
   X_1 - 2\Re (X_2) + X_3 &\ll \frac{(\log R)^2 (\log K)^3}{R} \sum_{d \leq R} \frac{\tau(d)^3}{\varphi(d)}\left( \frac{K}{d D} + R \right) \frac{K}{d}\\
&\ll(\log x)^5\left(\frac{K^2}{RD}+K\right).
\end{align*}
Since $R\leq x^{-\epsilon}NL$, $D\leq x^\epsilon$, and $K=NL$, we have $RD\leq K$, so the second term is absorbed by the first. Substituting the estimates for the three sums into~\eqref{eqn: dispersion cauchy schwarz}, and using $\widehat f(0)\ll M$ and $MK=x$, gives
\[
    \mcD^2\ll_\epsilon\frac{x^2(\log x)^6}{D}+x^{2-\delta}.
\]
In particular, under~\eqref{eqn: coprime support for triple convolution}, we have the stronger estimate
\begin{equation}\label{eqn: coprime convolution bound}
    \mcD\ll_\epsilon\frac{x(\log x)^3}{\sqrt D}+x^{1-\delta/2}.
\end{equation}

This proves the proposition in the special case~\eqref{eqn: coprime support for triple convolution}.

We now deal with the general case, where the coefficients need not satisfy~\eqref{eqn: coprime support for triple convolution}. Define
\begin{align*}
    \mcD(D):={}&\sum_{\substack{r\sim R\\(r,a_1a_2)=1}}
    \Biggl|\sum_{\substack{m\sim M\\n\sim N\\\ell\sim L}}
    \alpha_m\beta_n\gamma_\ell\\
    &\qquad\times\left(\mathbf1_{mn\ell\equiv a_1\overline{a_2}\pmod r}
    -\frac{\mathbf1_{(mn\ell,r)=1}}{\phi(r)}
    \omega_D(mn\ell\overline{a_1}a_2;r)\right)\Biggr|.
\end{align*}
Let $\delta_0>0$ be an exponent supplied by the special-case argument with $\epsilon/2$ in place of $\epsilon$. Choose
\[
    0<\delta\ll\eta\ll\min\{\epsilon,\delta_0\}.
\]
All constants here depend only on $\epsilon$; we decrease $\delta$ as necessary below.
We first reduce to $D\leq x^{\eta/4}$. If $D>x^{\eta/4}$, expanding $\omega_D-\omega_{x^{\eta/4}}$ into characters and applying the triangle inequality gives
\[
    |\mcD(D)-\mcD(x^{\eta/4})|
    \leq\sum_{r\sim R}\frac1{\phi(r)}
    \sum_{\substack{\chi\pmod r\\x^{\eta/4}<\operatorname{cond}(\chi)\leq D}}
    \left|\sum_{m\sim M}\alpha_m\chi(m)
    \sum_{n\sim N}\beta_n\chi(n)
    \sum_{\ell\sim L}\gamma_\ell\chi(\ell)\right|.
\]
We bound the $\ell$-sum trivially by $O(L)$. Write $r=dk$, where $d$ is the conductor of $\chi$, and let $\chi^*\pmod d$ be its inducing primitive character. Then
\[
    \chi(m)=\chi^*(m)\mathbf1_{(m,k)=1}.
\]
For each fixed $k$, the coefficients $\alpha_m\mathbf1_{(m,k)=1}$ and $\beta_n\mathbf1_{(n,k)=1}$ are $1$-bounded and independent of $d$. Using $\phi(dk)\geq\phi(d)\phi(k)$, Cauchy--Schwarz and the multiplicative large sieve~\cite[Lemme 6]{Drappeau} bound the contribution of a dyadic conductor block $E<d\leq2E$ by
\begin{align*}
    &\ll\frac{Lx^{o(1)}}E\sum_{k\leq2R/E}\frac1k
    \bigl(M(M+E^2)N(N+E^2)\bigr)^{1/2}\\
    &\ll\frac{x^{1+o(1)}}E
    \left(1+\frac E{\sqrt M}\right)
    \left(1+\frac E{\sqrt N}\right).
\end{align*}
Summing over $E$ therefore gives
\begin{equation}\label{eqn: change of conductor cutoff}
    |\mcD(D)-\mcD(x^{\eta/4})|
    \ll x^{1+o(1)}\left(x^{-\eta/4}+\frac1{\sqrt M}
    +\frac1{\sqrt N}+\frac D{\sqrt{MN}}\right).
\end{equation}
The hypotheses give $M\geq x^{1/3+5\epsilon}$, $N\geq x^\epsilon$, and $D\leq x^\epsilon$, so~\eqref{eqn: change of conductor cutoff} is $O(x^{1-c})$ for some $c=c(\epsilon)>0$. It therefore suffices to treat $D\leq x^{\eta/4}$, which we assume from now on.

We now show that the contribution to $\mcD(D)$ from terms with $(n\ell,a_1)>x^\eta$ is satisfactory. Define
\[
    b_v:=\sum_{\substack{mn\ell=v\\m\sim M,\ n\sim N,\ \ell\sim L\\(n\ell,a_1)>x^\eta}}
    \alpha_m\beta_n\gamma_\ell.
\]
By the triangle inequality, the contribution to $\mcD(D)$ from terms with $(n\ell,a_1)>x^\eta$ is at most $T_1+T_2$, where
\begin{align*}
    T_1&:=\sum_{\substack{r\sim R\\(r,a_1a_2)=1}}\;
    \sum_{\substack{v\\a_2v\equiv a_1\pmod r}}|b_v|,\\
    T_2&:=\sum_{\substack{r\sim R\\(r,a_1a_2)=1}}\frac1{\phi(r)}
    \left|\sum_{\substack{v\\(v,r)=1}}b_v\omega_D(v\overline{a_1}a_2;r)\right|.
\end{align*}
To bound $T_1$, note that the coefficients $b_v$ vanish unless $v\ll x$ and $(v,a_1)>x^\eta$, and satisfy $|b_v|\leq\tau_3(v)$. The divisor bound gives
\[
    \sum_v|b_v|\ll x^{1+o(1)}
    \sum_{\substack{g\mid|a_1|\\g>x^\eta}}\frac1g
    \ll x^{1-\eta+o(1)}.
\]
When $a_2v\neq a_1$, each modulus counted in $T_1$ divides $a_2v-a_1$; the possible term $a_2v=a_1$ contributes $O(Rx^{o(1)})$. Hence
\[
    T_1\ll x^{o(1)}(x^{1-\eta}+R).
\]

To bound $T_2$, expand $\omega_D$ into primitive characters of conductor $d$ and write $r=dk$, as in the cutoff comparison. Since a character modulo $d$ vanishes when $(v,d)>1$, the remaining restriction is $(v,k)=1$. Using $\phi(dk)\geq\phi(d)\phi(k)$ and the triangle inequality gives
\[
    T_2\ll x^{o(1)}\sum_{k\leq2R}\frac1k
    \sum_{d\leq D}\frac1d\sum_{\chi\pmod d}^{*}
    \left|\sum_{\substack{v\\(v,k)=1}}b_v\chi(v)\right|,
\]
where the star denotes primitive characters. For each fixed $k$, the coefficients $b_v\mathbf1_{(v,k)=1}$ are independent of $d$ and satisfy
\[
    \sum_{\substack{v\\(v,k)=1}}|b_v|^2
    \leq\sum_v|b_v|^2
    \ll x^{o(1)}\sum_v|b_v|
    \ll x^{1-\eta+o(1)}.
\]
For each fixed $k$, Cauchy--Schwarz followed by the multiplicative large sieve gives
\begin{align*}
    &\sum_{E<d\leq2E}\frac1d\sum_{\chi\pmod d}^{*}
    \left|\sum_{\substack{v\\(v,k)=1}}b_v\chi(v)\right|\\
    &\quad\leq\left(\sum_{E<d\leq2E}\frac{\phi(d)}{d^2}\right)^{1/2}
    \left(\sum_{E<d\leq2E}\sum_{\chi\pmod d}^{*}
    \left|\sum_{\substack{v\\(v,k)=1}}b_v\chi(v)\right|^2\right)^{1/2}\\
    &\quad\ll(E^2+x)^{1/2}
    \left(\sum_{\substack{v\\(v,k)=1}}|b_v|^2\right)^{1/2}
    \ll(E+\sqrt x)x^{(1-\eta)/2+o(1)}.
\end{align*}
Summing over dyadic conductor blocks and bounding the conductor-one contribution trivially gives
\begin{align*}
    T_2&\ll x^{(1-\eta)/2+o(1)}(D+\sqrt x)
    \sum_{k\leq2R}\frac1k\\
    &\ll x^{-\eta/2+o(1)}(x+D\sqrt x).
\end{align*}
Since $D\leq x^{\eta/4}$, we conclude that
\[
    T_1+T_2\ll x^{o(1)}(x^{1-\eta/2}+R).
\]
This is a power saving since $R\leq NL\leq x^{2/3-5\epsilon}$, so the contribution from $(n\ell,a_1)>x^\eta$ is satisfactory.

It remains to treat $(n\ell,a_1)\leq x^\eta$. Partition these terms by
\[
    d=(n,a_1),\qquad e=(\ell,a_1/d),\qquad
    g=de=(n\ell,a_1)\leq x^\eta,
\]
and write $n=dn'$, $\ell=e\ell'$. For fixed $d,e$, the new coefficients are
\[
    \beta^{(d)}_{n'}:=\beta_{dn'}\mathbf1_{(n',a_1/d)=1},\qquad
    \gamma^{(d,e)}_{\ell'}:=\gamma_{e\ell'}\mathbf1_{(\ell',a_1/g)=1}.
\]
They are $1$-bounded and satisfy~\eqref{eqn: coprime support for triple convolution} with $a_1/g$ in place of $a_1$. Since $g\mid|a_1|$ and $(r,a_1)=1$, division by $g$ in the residue class transforms this block into the preceding coprime case with
\[
    x'=x/g,\qquad N'=N/d,\qquad L'=L/e,\qquad
    a_1'=a_1/g,
\]
and with $M,R,a_2$ unchanged. Its original moduli form a subset of those coprime to $a_1'a_2$, so the bound for the latter set applies.

All parameter conditions hold with $x'$ and $\epsilon/2$ in place of $x$ and $\epsilon$. For example, $d,e\leq g\leq x^\eta$ gives
\[
    N'\geq x^{\epsilon-\eta}\geq(x')^{\epsilon/2},\qquad
    R\leq x^{-\epsilon}NL\leq(x')^{-\epsilon/2}N'L',
\]
provided $\eta\leq\epsilon/4$; substitution in the remaining inequalities in~\eqref{eqn: hypothesis for triple convolution} only improves their margins. Also,
\[
    |a_1'|\leq(x')^{1+\delta_0},\qquad
    |a_2|\leq(x')^{\delta_0},\qquad D\leq x^{\eta/4}\leq(x')^{\epsilon/2},
\]
by the choice of $\delta$ and $\eta$. We may therefore apply~\eqref{eqn: coprime convolution bound} to each block. The sum of its main bounds is controlled by
\[
    \sum_{\substack{d,e\geq1\\de\mid|a_1|}}\frac1{de}
    =\sum_{g\mid|a_1|}\frac{\tau(g)}g
    \leq\left(\frac{|a_1|}{\phi(|a_1|)}\right)^2
    \ll(\log\log x)^2.
\]
Thus the terms with $(n\ell,a_1)\leq x^\eta$ contribute
\[
    \ll_\epsilon\frac{x(\log x)^3(\log\log x)^2}{\sqrt D}+x^{1-c}
    \ll_\epsilon\frac{x(\log x)^4}{\sqrt D}+x^{1-c}
\]
for some $c=c(\epsilon)>0$, where the divisor bound absorbs the number of blocks in the power-saving error. Combining the three estimates and decreasing $\delta$ further so that $\delta\leq\eta/4$ and $\delta<2c$ yields, for the original range of $D$,
\[
    \mcD(D)\ll_\epsilon
    \frac{x(\log x)^4}{\sqrt{\min\{x^\delta,D\}}}.
\]
This proves the proposition without the auxiliary condition~\eqref{eqn: coprime support for triple convolution}.
\end{proof}


\begin{thebibliography}{14}
\bibitem{Assing-Blomer-Li}
E. Assing, V. Blomer, and J. Li. Uniform Titchmarsh divisor problems. \textit{Adv. Math.} \textbf{393} (2021), 108076. \href{https://doi.org/10.1016/j.aim.2021.108076}{doi:10.1016/j.aim.2021.108076}.
\bibitem{Baker-Harman}
R. C. Baker and G. Harman. Shifted primes without large prime factors. \textit{Acta Arith.} \textbf{83} (1998), 331--361. \href{https://doi.org/10.4064/aa-83-4-331-361}{doi:10.4064/aa-83-4-331-361}.
\bibitem{Balog}
A. Balog. On additive representation of integers. \textit{Acta Math. Hungar.} \textbf{54} (1989), 297--301.
\bibitem{Delabreteche-Tenenbaum}
R. de la Bretèche and G. Tenenbaum. Propriétés statistiques des entiers friables. \textit{Ramanujan J.} \textbf{9} (2005), 139--202.
\bibitem{Drappeau}
S. Drappeau. Théorèmes de type Fouvry-Iwaniec pour les entiers friables. \textit{Compos. Math.} \textbf{151} (2015), no.~5, 828--862.
\bibitem{Fouvry}
É. Fouvry. Autour du théorème de Bombieri-Vinogradov. II. \textit{Ann. Sci. École Norm. Sup.} \textbf{20} (1987), no.~4, 617--640.
\bibitem{HarperBV}
A. J. Harper. Bombieri--Vinogradov and Barban--Davenport--Halberstam type theorems for smooth numbers. Preprint (2012), arXiv:1208.5992v1 [math.NT]. \url{https://arxiv.org/abs/1208.5992v1}.
\bibitem{Hildebrand-Tenenbaum}
A. Hildebrand and G. Tenenbaum. On integers free of large prime factors. \textit{Trans. Amer. Math. Soc.} \textbf{296} (1986), no.~1, 265--290.
\bibitem{Kim-Sarnak}
H. H. Kim. Functoriality for the exterior square of $\text{GL}_4$ and the symmetric fourth of $\text{GL}_2$. \textit{J. Amer. Math. Soc.} \textbf{16} (2003), 139--183. With Appendix 1 by D. Ramakrishnan and Appendix 2 by H. H. Kim and P. Sarnak.
\bibitem{Lichtman}
J. D. Lichtman. Primes in arithmetic progressions to large moduli, and Goldbach beyond the square-root barrier. Preprint (2023), arXiv:2309.08522v1 [math.NT]. \url{https://arxiv.org/abs/2309.08522v1}.
\bibitem{Lichtman Shifted Primes}
J. D. Lichtman. Primes in arithmetic progressions to large moduli and shifted primes without large prime factors. Preprint (2022), arXiv:2211.09641v1 [math.NT]. \url{https://arxiv.org/abs/2211.09641v1}.
\bibitem{PascadiLargeSieve}
A. Pascadi. Large sieve inequalities for exceptional Maass forms and the greatest prime factor of $n^2+1$. \textit{Forum Math. Pi} \textbf{14} (2026), e8. \href{https://doi.org/10.1017/fmp.2026.10025}{doi:10.1017/fmp.2026.10025}.
\bibitem{PascadiLvlDist}
A. Pascadi. On the exponents of distribution of primes and smooth numbers. Preprint (2025), arXiv:2505.00653v2 [math.NT]. \url{https://arxiv.org/abs/2505.00653v2}.
\bibitem{PascadiSmoothBV}
A. Pascadi. Smooth numbers in arithmetic progressions to large moduli. \textit{Compos. Math.} \textbf{161} (2025), no.~8, 1923--1974.
\end{thebibliography}
\end{document}